\outer\def\give#1. {\medbreak
             \noindent{\bf#1. }}                     
\outer\def\section #1\par{\bigbreak\centerline{\S
     {\bf#1}}\nobreak\smallskip\noindent}
\def\({\left(}
\def\){\right)}

\def\sqr#1#2{{\vcenter{\hrule height.#2pt              
     \hbox{\vrule width.#2pt height#1pt\kern#1pt
     \vrule width.#2pt}
     \hrule height.#2pt}}}
\def\square{\mathchoice\sqr{5.5}4\sqr{5.0}4\sqr{4.8}3\sqr{4.8}3}
\def\qed{\hskip4pt plus1fill\ $\square$\par\medbreak}




\def\cB{{\cal B}}
\def\cC{{\cal C}}

\def\cG{{\cal G}}

\def\cI{{\cal I}}
\def\cJ{{\cal J}}

\def\cO{{\cal O}}

\def\cS{{\cal S}}
\def\cT{{\cal T}}

\def\cV{{\cal V}}
\def\cW{{\cal W}}


\def\CC{{\rm\kern.24em\vrule width.02em height1.4ex depth-.05ex\kern-.26em C}}
\def\RR{{\,\rm{\vrule width.02em height1.55ex depth-.07ex\kern-.3165em R}}}

\def\C{{\bf C}}
\def\cp1{{{\bf P}^1}}
\def\R{{\bf R}}

\def\Z{{\bf Z}}

\def\bar{\overline}              




\raggedbottom
\input epsf.sty
\magnification\magstep1
\centerline{\bf Real Polynomial Diffeomorphisms with Maximal Entropy:}
\centerline{\bf Tangencies}
\bigskip
\centerline {Eric Bedford* and John Smillie\footnote*{Research supported in
part by
the NSF.}}
\bigskip\bigskip
\section 0. Introduction

This paper deals with some questions about the dynamics of diffeomorphisms of
$\R^2$. A ``model family'' which has played a significant historical role
in dynamical systems and served as a focus for a great deal of research  is
the family introduced by H\'enon, which may be written as
$$f_{a,b}(x,y)=(a-by-x^2,x)\qquad b\ne0.$$
When $b\ne0$, $f_{a,b}$ is a diffeomorphism. When $b=0$ these maps
are essentially one dimensional, and the study of dynamics of
$f_{a,0}$ reduces to the study of the dynamics of quadratic maps
$$f_a(x)=a-x^2.$$ Like the H\'enon  diffeomorphisms of $\R^2$, the
quadratic maps of $\R$, have also played a central role in the field of
dynamical systems.

These two families of dynamical systems fit
together naturally, which raises the question of the extent to which the
dynamics is
similar. One difference is that our knowledge of these quadratic maps is
much more
thorough than our knowledge of these quadratic diffeomorphisms. Substantial
progress in the theory of quadratic maps has led to a rather complete
theoretical picture of their dynamics and an understanding of how
the dynamics varies with the parameter. Despite significant recent progress in
the theory of H\'enon diffeomorphisms, due to Benedicks and Carleson and many
others, there are still many phenomena that are not nearly so well
understood in
this two dimensional setting as they are for quadratic maps.

One phenomenon which illustrates the difference in the extent of our
knowledge in dimensions one and two is the dependence of the complexity of the
system on parameters. In one dimension the nature of this dependence is
understood,
and the answer is summarized by the principle of monotonicity.  Loosely stated,
this is the assertion that the complexity of $f_a$ does not decrease as the
parameter $a$ increases. The notion of complexity used here can be made precise
either in terms of counting periodic points or in terms of entropy. The paper
[KKH] shows that monotonicity
is a much
more complicated phenomenon for diffeomorphisms. In this paper we will
focus on one
end of the complexity spectrum, the diffeomorphisms of maximal entropy, and
we will
show to what extent the dynamics in the two dimensional case are similar to the
dynamics in the one dimensional case. In the case of quadratic maps complex
techniques proved to be an important tool for developing the theory. In
this paper
we apply complex techniques to study quadratic (and higher degree)
diffeomorphisms.

Topological entropy is a measure of dynamical complexity that can be
defined
either for maps or diffeomorphisms. By Friedland and Milnor [FM] the
topological
entropy of H\'enon diffeomorphisms satisfies: $0\le h_{top}(f_{a,b})\le\log
2$.
We will say that $f$ has maximal entropy if the topological entropy is equal to
$\log 2$. The notion of maximal entropy makes sense for polynomial maps of $\R$
as well as polynomial diffeomorphisms of $\R^2$ of degree greater than two. In
either of these cases we say that $f$ has maximal entropy if
$h_{top}(f)=\log(d)$ where $d$ is the algebraic degree of $f$ and $d\ge2$. We
will see that this condition is equivalent to the assumption that $f^n$ has
$d^n$ (real) fixed points for all $n$.

The quadratic maps $f_a$ of maximal entropy are those with $a\ge2$.
These maps are hyperbolic (that is to say expanding) for $a>2$, whereas the
map $f_2$, the example of Ulam and von Neumann, is not hyperbolic.
Examples of maps of maximal entropy in the H\'enon family were given by Devaney
and Nitecki [DN] (see also [HO] and [O]), who showed that for certain parameter
values
$f_{a,b}$ is hyperbolic and a model of the Smale horseshoe. Examples of maximal
entropy polynomial diffeomorphisms of degree $d\ge2$ are given by
the $d$-fold horseshoe mappings of Friedland and Milnor (see [FM, Lemma 5.1]).

We will see that all polynomial diffeomorphisms of maximal entropy (whether or
not they are hyperbolic) exhibit a certain form of expansion.
Hyperbolic diffeomorphisms have uniform expansion and contraction which implies
uniform expansion and contraction for periodic orbits. To be precise, let
$p$ be a point of period $n$ for a diffeomorphism $f$. We say that $p$ is a
saddle point if $Df^n(p)$ has eigenvalues $\lambda^{s/u}$ with
$|\lambda^s|<1<|\lambda^u|$. If $f$ is hyperbolic then for some $\kappa>1$
independent of $p$ we have $|\lambda^u|\ge\kappa^n$ and
$|\lambda^s|\le\kappa^n$. On the other hand it is not true that uniform
expansion/contraction for periodic points implies hyperbolicity. A one
dimensional
example of a map with expansion at periodic points which is not hyperbolic
is given
by the Ulam-von Neumann map. This map is not expanding because the critical
point, $0$, is contained in the non-wandering set, $[-2,2]$ (which is also
equal to
$K$). The map satisfies the above inequalities with
$\kappa=2$. In fact for every periodic point of period $n$ except the fixed
point
$p=-2$ we have $|Df^n(p)|=2^n$. At $p=-2$ we have $n=1$ yet $|Df^n(p)|=4$.

\proclaim Theorem 1. If $f$ is a maximal entropy polynomial diffeomorphism,
then
\item{(1)} Every periodic point is a saddle point.
\item{(2)} Let $p$ be a periodic point of period $n$. Then $|\lambda^s(p)|<
1/d^{n}$, and
$|\lambda^u(p)|> d^n$.
\item{(3)} The set of fixed points of $f^n$ has exactly $d^n$ elements.

  Let $K$ be the set of points in $\R^2$ with
bounded orbits. In Theorem 5.2 (below) we show that $K$ is a Cantor set for
every
maximal entropy diffeomorphism.  By [BS8, Proposition 4.7] this yields
the strictness
of the inequalities in (2). Note that the situation for maps of maximal entropy
in one variable is different. In the case of the Ulam-von Neuman map $K$ is
a connected interval and the strict inequalities do not hold.

We note that by [BLS], condition (3) implies that $f$ has maximal entropy.
Thus
we see that condition (3) provides a way to characterize the class of maximal
entropy diffeomorphisms which makes no explicit reference to entropy. As
was noted
above, we can define the set of maximal entropy diffeomorphisms using either
notion of complexity: they are the polynomial diffeomorphisms for which entropy
is as large as possible, or equivalently those having as many periodic
points as
possible.

For the Ulam-von Neumann map the fixed point $p=-2$ which is the left-hand
endpoint
of $K$ is distinguished as was noted above. This distinction has an analog in
dimension two. Let $p$ be a saddle point. Let
$W^{s/u}(p)$ denote the stable and unstable manifolds of $p$. These sets
are analytic
curves. We say a periodic point $p$ is $s/u$ one-sided if only one component of
$W^{s/u}-\{p\}$ meets $K$.
For one-sided periodic points the estimates of Theorem 1 (2) can be
improved. If $p$ is $s$ one-sided, then
$|\lambda^s(p)|<1/d^{2n}$; and if
$p$ is
$u$ one-sided, then $|\lambda^u(p)|>d^{2n}$.

The set of parameter values corresponding to diffeomorphisms of maximal
entropy is closed, while the set of parameter values corresponding to
hyperbolic
diffeomorphisms is open.  It follows that not all maximal entropy
diffeomorphisms
are hyperbolic.  We now address the question: which properties of hyperbolicity
fail in these cases.

\proclaim Theorem 2. If $f$ has maximal entropy, but $K$ is not a hyperbolic
set for $f$, then
\item{(1)} There are periodic points $p$ and $q$ in $K$ (not necessarily
distinct) so that $W^u(p)$ intersects $W^s(q)$ tangentially with order 2
contact.
\item{(2)} $p$ is $s$ one-sided, and $q$ is $u$ one-sided.
\item{(3)} The restriction of $f$ to $K$ is not expansive.

Condition (1) is incompatible with $K$ being a hyperbolic set. Thus this
theorem
describes a specific way in which hyperbolicity fails. Condition (3), which is
proved in [BS8, Corollary 8.6], asserts that for any
$\epsilon>0$ there are points $x$ and $y$ in $K$ such that for all
$n\in\Z$, $d(f^n(x),f^n(y))\le\epsilon$. Condition (3) is a topological
property which is not compatible with hyperbolicity.  We conclude
that when $f$ is not hyperbolic it is not even {\it topologically} conjugate to
any hyperbolic diffeomorphism.

The proofs of the stated
theorems owe much to the theory of quasi-hyperbolicity developed in [BS8].
In [BS8] we show that
maximal entropy diffeomorphisms are quasi-hyperbolic.
We also define a singular set $\cC$ for any quasi-hyperbolic
diffeomorphism.  Much of the work of this paper is devoted to
showing that
in the maximal entropy case
$\cC$ is finite and consists of one-sided periodic points. Further
analysis allows us to show that they these periodic points have period
either 1 or
2. In the case of quadratic mappings we can say exactly which points are
one-sided.

We say that a saddle point is non-flipping if $\lambda^u$ and $\lambda^s$ are
both positive.

\proclaim Theorem 3.  Let $f_{a,b}$ be a quadratic mapping with maximal
entropy.  If $f_{a,b}$ preserves orientation, then the unique non-flipping
fixed
point of $f$ is doubly one-sided.  If $f$ reverses orientation, then one of its
fixed points is stably one-sided, and the other is unstably one-sided.  There
are no other one-sided points in either case.

 We can use our results to describe how hyperbolicity is lost on the
boundary of
the horseshoe region for H\'enon diffeomorphisms. We focus on the orientation
preserving case here, but our results allow us to treat the orientation
preserving case as well. The parameter space for orientation preserving H\'enon
diffeomorphisms is the set $\{(a,b): b>0\}$. Let us define the horseshoe
region
to be the largest connected open set containing the Devaney-Nitecki horseshoes
and consisting of hyperbolic diffeomorphisms.  Let
$f=f_{a_0,b_0}$ be a point on the boundary of the horeshoe region. It follows
from the continuity of entropy that $f$ has maximal entropy. Theorem 1 tells us
that $f$ has the same number of periodic points as the horseshoes and that they
are all saddles. In particular no bifurcations of periodic points occur at
${a_0,b_0}$. Let
$p_0$ be the unique non-flipping fixed point for $f$. It follows from Theorem 2
that the stable and unstable manifolds of $p_0$ have a quadratic homoclinic
tangency.

\epsfxsize5.5in
\centerline{\epsfbox{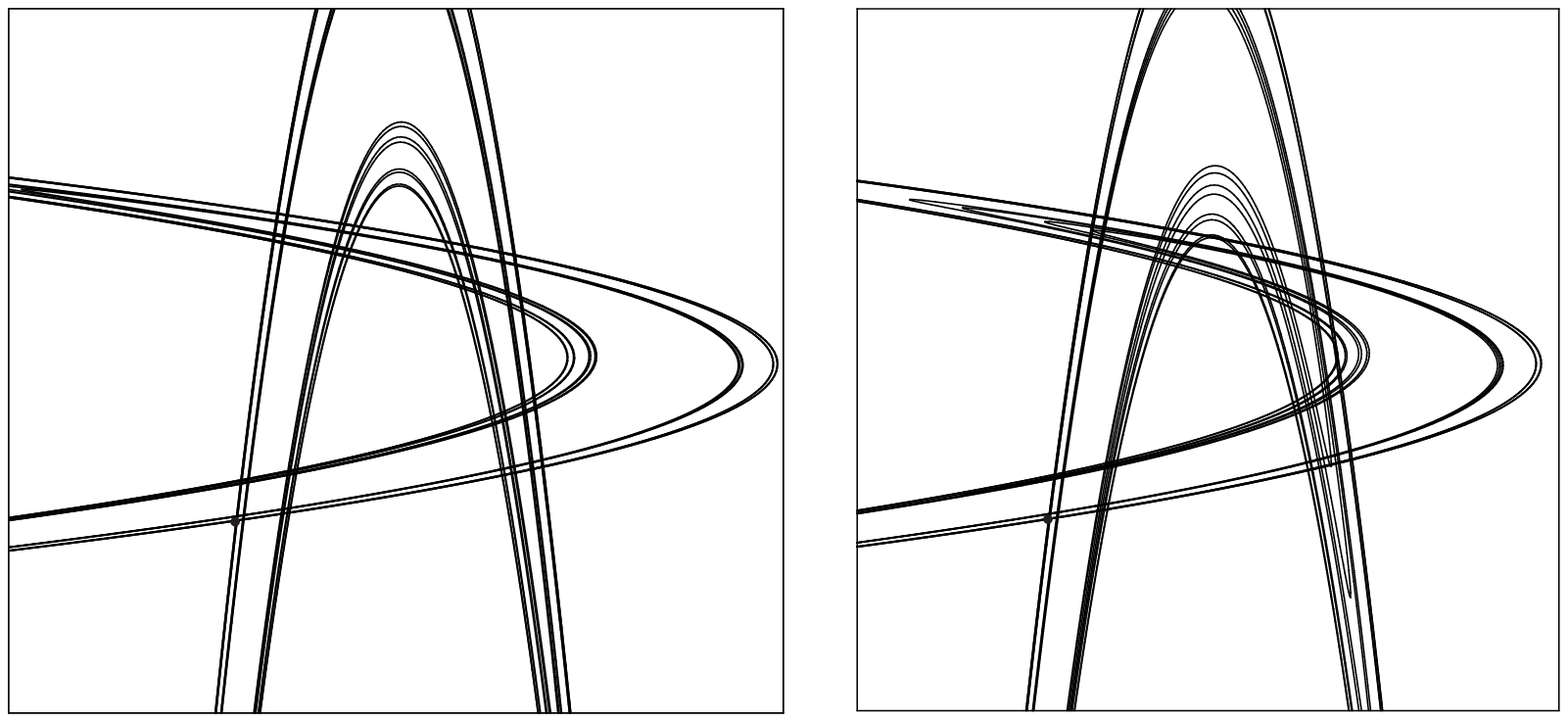}}

\centerline{Figure 0.1}

Figure 0.1 shows computer-generated pictures of mappings $f_{a,b}$ with
$a=6.0$, $b=0.8$ on the left and $a=4.64339843$, $b=0.8$ on the
right.\footnote*{We thank Vladimir Veselov for using a computer program
that he wrote to generate this second set of parameter values for us.} The
curves
pictured are the stable/unstable manifolds of the non-flipping saddle point
$p_0$, which is the point marked by a disk in each picture at the lower
leftmost
point of intersection of the stable and unstable manifolds.  The manifolds
themselves are connected; the apparent disconnectedness is a result of clipping
the picture to a viewbox. There are no tangential intersections evident on the
left, while there appears to be a tangency on the right.  This is consistent
with the analysis above.

\section 1. Background

Despite the fact that we study real polynomial diffeomorphisms, the proofs
of the results of this paper depend on the theory of complex polynomial
diffeomorphisms. In particular the theory of quasi-hyperbolicity which lies at
the heart of much of what we do is a theory of complex polynomial
diffeomorphisms. The notation we use in the paper is chosen to simplify the
discussion of complex polynomial diffeomorphisms.  A polynomial diffeomorphisms
of $\C^2$ will be denoted by
$f_\C$, or simply $f$, when no confusion will result.  Let $\tau(x,y)=(\bar
x,\bar y)$  denote complex conjugation in $\C^2$.  The fixed point set of
complex
conjugation in $\C^2$ is exactly $\R^2$.  We say that $f$ is real when
$f:\C^2\to\C^2$ has real coefficients, or equivalently, when $f$ commutes with
$\tau$.  When $f$ is real we write $f_\R$ for the restriction of $f$ to
$\R^2$.

  Let us consider mappings of the form $f=f_1\circ\cdots\circ f_m$,
where
$$f_j(x,y)=(y,p_j(y)-a_j x),\eqno(1.1)$$
$p_j$ is a polynomial of degree $d_j$. If we set $d=d_1\ldots d_m$, then it
is easily seen that if $f$ has the form 1.1 then the degree of $f$ is $d$. The
degree of $f^{-1}$ is also $d$ and, since $h(f_\R)=h(f^{-1}_\R)$ it follows
that
$f$ has maximal entropy if and only if $f^{-1}$ does.

\proclaim Proposition 1.1. If a real polynomial diffeomorphism $f$ has
maximal entropy, then it is conjugate via a real polynomial diffeomorphism to a
real polynomial diffeomorphism of the same degree in the form (1.1).

\give Proof. According to [FM] a polynomial diffeomorphism $f_\R$ of $\R^2$ is
conjugate via a polynomial diffeomorphism, $g$, to an elementary diffeomorphism
or a diffeomorphism in form (1.1). Since $f_\R$ has
positive entropy it is not conjugate to an elementary diffeomorphism. In
[FM] it
is also show that a diffeomorphism in the form (1.1) has minimal entropy among
all elements in its conjugacy class so
$deg(g_\R)\le deg(f_\R)$. Since entropy is a conjugacy invariant we have:
$$\log deg(g_\R)\le\log deg(f_\R)=h(f_\R)=h(g_\R).$$ Again by [FM],
$h(g_\R)\le\log
deg(g_\R)$ so we conclude that the inequalities are equalities and that
$deg(g_\R)=
deg(f_\R)$.

Thus we may assume that we are dealing with maximal entropy polynomial
diffeomorphisms written in  form (1.1). The mapping $f_{a,b}$ in the
Introduction is not in the form (1.1), but the linear map $L(x,y)=(-y,-x)$
conjugates $f_{a,b}$ to
$$(x,y)\mapsto(y,y^2-a-bx).$$
In Sections 1 through 4, we are dealing with polynomial diffeomorphisms of all
degrees, and we will assume that the are given in the form (1.1).

We recall some standard notation for general polynomial diffeomorphisms of
$\C^2$.  The set of points in $\C^2$ with bounded forward orbits is denoted by
$K^+$. The set of points with bounded backward orbits is denoted by $K^-$. The
sets $J^\pm$ are defined to be the boundaries of $K^\pm$. The set $J$ is
$J^+\cap
J^-$ and the
set $K$ is $K^+\cap K^-$. Let $\cS$ denote the set of saddle points of $f$.
For a
general polynomial diffeomorphism of $\C^2$ the closure of $\cS$ is denoted by
$J^*$. For real polynomial diffeomorphism of $\C^2$ each of these $f$-invariant
sets is
also invariant under $\tau$. For a real maximal entropy mapping it is proved in
[BLS] that $J^*=J=K$ and furthermore that this set is real; that is
$K\subset\R^2$.

For $p\in\cS$, there is a holomorphic immersion
$\psi^u_p:\C\to\C^2$ such that $\psi^u_p(0)=p$ and $\psi^u_p(\C)=W^u(p)$. The
immersion $\psi^u_p$ is well defined up to multiplication by a non-zero complex
scalar. By using a certain potential function we can choose  distinguished
parametrizations. Define $G^+$ by the formula
$$G^+(x,y)=\lim_{n\to\infty}{1\over d^n}\log^+|f^n(x,y)|.$$ Changing parameter
in the domain via a change of coordinates $\zeta'=\alpha\zeta$,
$\alpha\ne0$, we may assume that $\psi^u_p$ satisfies
$$\max_{|\zeta|\le1}G^+\circ\psi^u_p(\zeta)=1.$$ With this normalization,
$\psi^u_p$ is uniquely determined modulo rotation; that is, all such mappings
are of the form $\zeta\mapsto\psi_p^u(e^{i\theta}\zeta)$.

When the diffeomorphism $f$ is real and $p\in\R^2$ we may choose the
parametrization of $W^u_p$ so that it is real, which is to say that
$\psi=\psi^u_p$ satisfies $\psi(\bar\zeta)=\tau\circ\psi(\zeta)$. In this case
the set $\psi^{-1}(K)=\psi^{-1}(K^+)$ is symmetric with respect to the real
axis
in $\C$ and the parametrization is well defined up to multiplication
by $\pm1$. In the real case $\psi(\R)\subset\R^2$, and the set $\psi(\R)$ is
equal to the unstable manifold of $p$ with respect to the map $f_\R$.

When $f$ is real and has maximal entropy more is true. In this case every
periodic point is contained in $\R^2$. Let $\psi$ be a real parametrization.
Since $\psi$ is injective,  the inverse image of the fixed point set in $\C^2$
is contained in the fixed point set in $\C$.  Thus $\psi^{-1}(\R^2)=\R$, and
$\psi^{-1}(K)\subset\R$.  If $p$ is a $u$ one-sided periodic point then $K$
meets only one component of $W^u(p,\R)$ so $\psi^{-1}(K)$ is contained in
one of
the rays $\{\zeta\in\R:\zeta\ge0\}$ or $\{\zeta\in\R:\zeta\le0\}$.

We define the set of all such unstable parametrizations as
$\psi^u_\cS:=\{\psi_p^u:p\in\cS\}$.  For $\psi\in\psi^u_p$ there exist
$\lambda=\lambda^u_{p}\in\R$ and $\tilde f\psi\in\psi^u_{fp}$ such that
$$(\tilde f\psi)(\zeta)=f(\psi(\lambda^{-1}\zeta))\eqno(1.2)$$
for $\zeta\in\C$.

A  consequence of the fact that $\psi^{-1}(K)\subset\R$ [BS8, Theorem 3.6] is
that
$$|\lambda_p|\ge d.\eqno(1.3)$$
Furthermore if $p$ is u one-sided then
$$|\lambda_p|\ge d^2.$$
The condition that $|\lambda_p|$ is bounded below by a constant greater
than one
is one of several equivalent conditions that can serve as definitions of the
property of quasi-expansion defined in [BS8]. Thus, as in [BS8], we see
that $f$
and $f^{-1}$ are quasi-expanding. A consequence of quasi-expansion is that
$\psi^u_\cS$ is a normal family (see [BS8, Theorem 1.4]).  In this case we
define
$\Psi^u$ to be the set of normal (uniform on compact subsets of
$\C$) limits of elements of $\psi^u_\cS$.  Let
$\Psi^u_p:=\{\psi\in\Psi^u:\psi(0)=p\}$. It is a further consequence of
quasi-expansion that $\Psi^u$ contains no constant mappings.

For $p\in J$, the mappings in $\Psi^{u}_p$ have a common image which we
denote by
$V^u(p)$ ([BS8, Lemma 2.6]). Let $W^u(p)$ denote the ``stable
set'' of
$p$. This
consists of
$q$ such that $$\lim_{n\to+\infty}dist(f^{-n}p,f^{-n}q)=0.$$
It is proved in  [BS8, Proposition 1.4] that $V^u(p)\subset W^u(p)$. It follows
that $V^u(p)\subset K^-$. In many cases the stable set is actually one
dimensional complex manifold. When this is the case it follows that
$V^u(p)=W^u(p)$.

Let $V^u_\epsilon(p)$ denote the component of $V^u(p)\cap B(p,\epsilon)$ which
contains $p$. For $\epsilon$ sufficiently small $V^u_\epsilon(p)$ is a properly
embedded variety in $B(p,\epsilon)$. Let $E^{u}_p$ denote the tangent space
to this
variety at $p$. It may be that the variety $V^u_\epsilon(p)$ is
singular at $p$. In this case we define the tangent cone to be the set of
limits of
secants.

For $\psi\in\Psi^{u}$ we say that ${\rm Ord}(\psi)=1$ if
$\psi'(0)\ne0$; and if
$k>1$, we say ${\rm Ord}(\psi)=k$ if $\psi'(0)=\cdots=\psi^{(k-1)}(0)=0$,
$\psi^{(k)}(0)\ne0$. Since $\Psi^u$ contains no constant mappings, ${\rm
Ord}(\psi)$
is finite for each $\psi$. If
$\psi\in\Psi^{s/u}$, and if
${\rm Ord}(\psi)=k$, then there are $a_j\in\C^2$ for $k\le j<\infty$ such that
$$\psi(\zeta)=p+a_k\zeta^k+a_{k+1}\zeta^{k+1}+\dots$$

  It is easy to show that the tangent cone $E^{u}_p$ to the
variety $V^u_\epsilon(p)$ is actually the complex subspace of the tangent space
$T_p\C^2$ spanned by $a_k$. One consequence of this is that
the span of the $a_k$
term depends only on $p$ and not on the particular mapping in $\Psi_p^{u}$.
(It is
possible however that different parametrizations give different values for
$k$.)
A second consequence is that even when the variety $V^u_\epsilon(p)$ is
singular
the tangent cone is actually a complex line and, in what follows, we will refer
to $E^{u}_p$ as the tangent space. The mapping
$\psi\mapsto{\rm Ord}(\psi)$ is an upper semicontinuous function on
$\Psi^{u}$.  For $p\in J$, we set
$\tau^{u}(p)=\max\{{\rm Ord}(\psi):\psi\in\Psi^{u}_p\}$.
The reality of $\psi$ is equivalent to the
condition that $a_j\in\R^2$.

Since $f^{-1}$ is also
quasi-expanding, we may repeat the definitions above with $f$ replaced by
$f^{-1}$ and unstable manifolds replaced by stable manifolds; and in this
case we
replace the superscript $u$ by $s$.   We set
$$J_{j,k}=\{p\in J:\tau^s(p)=j,\tau^u(p)=k\}.$$
We define
$$\lambda^{s/u}(p,n)=\lambda^{s/u}_{p}\cdots\lambda^{s/u}_{f^{n-1}p}.$$
Iterating the mapping $\tilde f$ defined above, we have mappings $\tilde
f^n:\Psi^{s/u}_p\to\Psi^{s/u}_{f^np}$ defined by
$$\tilde f^n(\psi^{s/u}(\zeta))=f^n\circ
\psi^{s/u}(\lambda^{s/u}(p,n)^{-1}\zeta).\eqno(1.5)$$
By (1.3),
$$|\lambda^s(p,n)|\le d^{-n}, \ \ \ \  |\lambda^u(p,n)|\ge d^n.\eqno(1.6)$$

We will give here the proof of Theorem 3.
Since $f$ and $f^{-1}$ are quasi-expanding it follows that every periodic
point in
$J^*$ is a saddle. Since every periodic point is contained in $K$ and
$K=J^*$ it
follows that every periodic point is a saddle. According to [FM] the number
of fixed
points of $f_\C^n$ counted with multiplicity is $d^n$. Since all
periodic points are saddles they all have multiplicity one (multiplicity is
computed with respect to $\C^2$ rather than $\R^2$). Thus the set of fixed
points
of
$f^n$ has cardinality
$d^n$. Since
$K\subset\R^2$ all of
these points are real.

\section 2.  The Maximal Entropy Condition and its Consequences

Let us return to our discussion of the maximal entropy condition. The
argument that
$\psi^{-1}(\R^2)=\R$ depended on the injectivity of $\psi$. Even though
elements
of
$\Psi^u$ are obtained by taking limits of elements of $\psi^u_\cS$ it does not
follow that
$\psi\in\Psi^u$ is injective.  In fact it need not be the case that
$\psi^{-1}(\R^2)\subset\R$, but the following Proposition shows that a related
condition still holds.

\proclaim Proposition 2.1.  For $\psi\in\Psi^u$, $\psi^{-1}(K)\subset\R$.

\give Proof. The image of $\psi$ is contained in $K^-$, it follows that
$\psi^{-1}(K^+)=\psi^{-1}(K)$ for $\psi\in\psi^u_\cS$.  Since
$G^+$ is harmonic on $\C^2-K^+$, it follows that $G^+\circ\psi$ is harmonic on
$\C-\R\subset\C-\psi^{-1}K$.  By Harnack's principle, $G^+\circ\psi$ is
harmonic on
$\C-\R$ for any limit function $\psi\in\Psi^u$.  If $G^+\circ\psi$ is zero
at some
point $\zeta\in\C-\R$ with, say, $\Im(\zeta)>0$, then it is zero on the
upper half
plane by the minimum principle.  By the invariance under complex conjugation,
it is zero everywhere.  But this means that $\psi(\C)\subset\{G^+=0\}=K^+$.  By
(1.4), this means that $\psi(\C)\subset K\subset\R^2$.  Since $K$ is bounded,
$\psi$ must be constant.  But this is a contradiction because $\Psi^u$ contains
no constant mappings.
\qed

Our next objective is to find a bound on ${\rm Ord}(\psi)$ for $\psi\in\Psi^u$.
Set $m^{u}=\max_J\tau^{u}$ and consider the maximal index $j$ so that
$J_{j,m^u}$ is
non-empty. Thus $J_{j,m^u}$ is a maximal index pair in the language of
[BS8]. By
[BS8, Proposition 5.2],
$J_{j,m^u}$ is a hyperbolic set with stable/unstable subspaces given by
$E^{s/u}_p$.

The notion of a homogeneous parametrization was defined in [BS8, \S6].
A homogeneous parametrization of order $m$, $\psi:\C\to\C^2$, is one that can
be written as $\psi(\zeta)=\phi(a\zeta^m)$ for some $a\in\C-\{0\}$ and some
non-singular $\phi:\C\to\C^2$. It follows from [BS8, Lemma 6.5] that for every
$p$ in a maximal index pair such as $J_{j,m^u}$ there is homogeneous
parametrization in $\Psi^u_p$ with order $m^u$.

\proclaim Proposition 2.2.  Suppose that $\psi\in\Psi^u$, is a homogeneous
parametrization of order $m$.  Then it follows that $m\le 2.$

\give Proof.  By Proposition 2.1, $\psi^{-1}(J)\subset\R$.  And from the
condition $\psi(\zeta)=\phi(\zeta^m)$ it follows that $\psi^{-1}(J)$ is
invariant
under rotation by $m$-th roots of unity.  Now $\psi^{-1}(J)$ is non-empty
(containing 0) and a non-polar subset of $\C$, since it is the zero set of the
continuous, subharmonic function $G^+\circ\psi$. Since a polar set contains no
isolated points it follows that
$\psi^{-1}(J)$ contains a point $\zeta_0\ne0$.  Since the rotations of
$\zeta_0$ by the $m$-th roots of unity must lie in $\R$, it follows that
$m\le2$. \qed

\proclaim Corollary 2.4.
$J=J_{1,1}\cup J_{2,1}\cup J_{1,2}\cup J_{2,2}.$

There are three possibilities to consider.

\item{(1)} $J_{2,1}\cup J_{1,2}\cup J_{2,2}$ is empty. In this case it
follows from [BS8] that $f$ is hyperbolic.

\item{(2)} $J_{2,2}$ is empty and $J_{2,1}\cup J_{1,2}$ is non-empty. In
this case
$J_{2,1}$ and $J_{1,2}$ are maximal index pairs and are both hyperbolic sets.

\item{(3)} $J_{2,2}$ is non-empty but $J_{2,1}\cup J_{1,2}$ is empty. In
this case
$J_{2,2}$ is a maximal set and is hyperbolic.

\item{(4)} $J_{2,2}$ is non-empty and  $J_{2,1}\cup J_{1,2}$ is non-empty.
This is the only case in which we do not know a priori that points in
$J_{2,1}\cup J_{1,2}$ are regular.

\proclaim Proposition 2.5.  For $p\in J$, let $\psi\in\Psi^{s/u}_p$ be
given.  Then
$\zeta\mapsto\psi(\zeta)$ is at most two-to-one.  If $\psi$ is two-to-one, then
it has one critical point, which must be real.

\give Proof.  This  follows from Proposition 2.2 and [BS8, Lemma 4.6].
\qed

\proclaim Proposition 2.6.  Let $p$ be in $J_{*,2}$, and let
$V^u_\epsilon(p)$ be regular.
If $\psi\in\Psi^u_p$ has order 2, there is an embedding $\phi:\C\to\C^2$
such that
$\psi(\zeta)=\phi(\zeta^2)$.

\give Proof.  By Proposition 2.5, $\psi$ has at most one critical point,
which must
be $\zeta=0$.  Thus all points of $\psi(\C)-\{p\}$ are regular.  Since
$V^u_\epsilon(p)$ is
regular, it follows that $\psi(\C)$ is regular, so there is an embedding
$\phi:\C\to\C^2$ with $\phi(\C)=\psi(\C)$.  By Proposition 2.3 $\tau^u\le
2$, so
$J_{*,2}$, being a set of maximal order, is compact.  Thus $\alpha(p)\subset
J_{*,2}$, so the result follows from [BS8, Proposition 4.4]. \qed

If
$\psi\in\Psi^u_p$ is one-to-one, then $\psi(\C)\cap\R^2=\psi(\R)$.  (For if
there is
a point $\zeta\in\C-\R$ with $\psi(\zeta)\in\R^2$, then we would also have
$\psi(\bar\zeta)\in\R^2$.  But $\zeta\ne\bar\zeta$, contradicting the
assumption
that $\psi$ is one-to-one.)  If $\psi$ is 2-to-1, then $\psi$ has a
critical point
$t_0\in\R$.  Let us suppose that $\psi$ has a quadratic singularity at
$\zeta=0$,
i.e. $\psi(\zeta)=p+a_2\zeta^2+O(|\zeta|^3)$.  If $\psi(\C)\cap\R^2$ is a
smooth
curve, then $p$ divides this curve into two pieces: in Figure 2.1 the image of
$\R$ under $\psi$ is drawn dark, and the image of $i\R$ is shaded.  By
Proposition 2.1, the shaded region is disjoint from $J$.
\epsfxsize2.1in
\bigskip
\centerline{\epsfbox{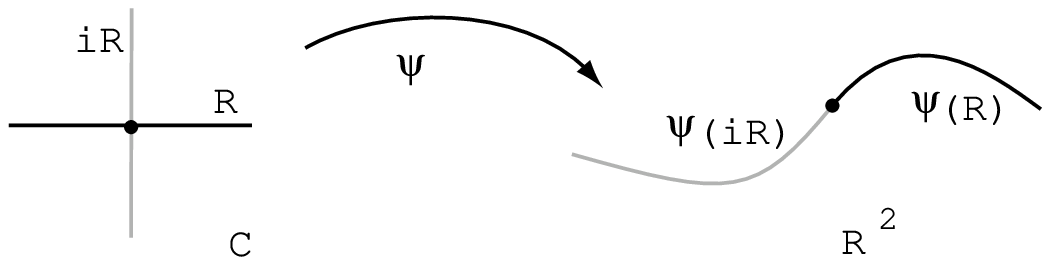}}
\centerline{Figure 2.1}

   Recall that the tangent space to $V^{s/u}_\epsilon(p)$ at $p$
is $E^{s/u}_p$.  We say that $V^u_\epsilon(p)$ and $V^s_\epsilon(p)$
intersect tangentially at $p$ if
$E^s_p=E^u_p$.  We recall that $\alpha(p)$, the $\alpha$-limit set of $p$, is
the set of limit points of $\{f^{-n}p: n\ge0\}$, and the $\omega$-limit
set, $\omega(p)$, is the set of limit points of $\{f^np:n\ge0\}$.
Compactness of
$J$ implies that $\alpha(p)$ and $\omega(p)$ are non-empty.  The following are
consequences of Theorem 7.3 of [BS8].

\proclaim Theorem 2.7.  Suppose the varieties $V^u_\epsilon(p)$ and
$V^s_\epsilon(p)$ intersect
tangentially at $p\in J$ (i.e.\ suppose $E^s_p=E^u_p$).  Then the $\alpha$- and
$\omega$-limit sets satisfy $\alpha(p)\subset J_{2,*}$ and
$\omega(p)\subset J_{*,2}$.  Further, $p$ belongs to $J_{1,1}$, and the
varieties of $V_p^{s/u}$ are regular at $p$.

\proclaim Theorem 2.8.  If $V^s_\epsilon(p)$ and $V^u_\epsilon(p)$ are
tangent at
$p\in J$, then the tangency is at most second order, i.e.\
$V^s_\epsilon(p)$ and
$V^u_\epsilon(p)$ have different curvatures at $p$.

\section 3.  Finiteness of Singular Points

Let us consider a point $p\in J$ where the varieties $V^s_\epsilon(p)$ and
$V^u_\epsilon(p)$ are
nonsingular and intersect transversally.  We may perform a real, affine
change of
coordinates so that in the new coordinate $(x,y)$ we have $p=(0,0)$,
$V^u_\epsilon(p)$ is
tangent to the
$x$-axis at
$p$, and
$V^s_\epsilon(p)$ is tangent to the $y$-axis at $p$.  let $\pi_s(x,y)=y$
and $\pi_u(x,y)=x$.
For $\epsilon>0$ let $\Delta(\epsilon)=\{\zeta\in\C:|\zeta|<\epsilon\}$.  For
$q\in\Delta^2(\epsilon)\cap J$ let $V^{s/u}(q,\epsilon)$ denote the
connected component of
$V^{s/u}_q\cap\pi_{s/u}^{-1}\Delta(\epsilon)$ containing $q$.  For $\epsilon>0$
small we have
$$\pi_s:V^u(p,\epsilon)\subset\Delta(\epsilon/2),\
\ \pi_u:V^s(p,\epsilon)\subset\Delta(\epsilon/2),\eqno(3.1)$$ and
$$\pi_{s/u}:V^{s/u}(p,\epsilon)\to\Delta(\epsilon){\rm\ are\ proper\ maps\ of\
degree\ 1.}\eqno(3.2)$$
By [BS8, Lemmas 2.1 and 2.2] the  varieties $V^{s/u}_\epsilon(q)$ depend
continuously
on
$q$.  Thus for $\delta>0$ small, (3.1) will hold for the varieties at $q$ if
$q\in\Delta^2(\delta)\cap J$, and  the projections
$\pi_{s/u}:V^{s/u}(q,\epsilon)\to\Delta(\epsilon)$ will be proper.

By a {\it regular box } $\cB$ we will refer to the ensemble of affine
coordinate system, projections $\pi_{s/u}$, and sets
$\Delta^2(\epsilon)$, $\Delta^2(\delta)$.  To begin with, we require that
$\epsilon$
and
$\delta$ are small enough that (3.1) and (3.2) hold.  Then we will shrink
$\epsilon$ and  $\delta$ progressively so that the subsequent results from
Lemmas
3.1 through 3.5 hold.  Let us define $\cV^s(\cB)$ as the set of
varieties
$V^s(q,\epsilon)$ for $q\in\bar\Delta^2(\delta)\cap J$.  Further, we define
$\cV^s_j$ as the set of varieties $V^s\in\cV^s(\cB)$ such that the projection
$\pi_s|_{V^s}:V^s\to\Delta(\epsilon)$ has mapping degree $j$.  In a similar
way, we
define $\cV^u(\cB)$ and $\cV^u_j(\cB)$.  It is evident that elements of
$\cV_1^{s/u}(\cB)$ are represented as graphs of analytic functions, and so
$\cV_1^{s/u}(\cB)$ is a compact family of varieties.

\proclaim Lemma 3.1.  If
$V^s\in\cV^s_j$, $V^u\in\cV^u_k$, then the intersection $V^s\cap V^u$
consists of $jk$ points (counted with ``intersection'' multiplicity). If
$\epsilon$ and $\delta$ are sufficiently small, then
$\cV^s=\cV^s_1\cup\cV^s_2$.

\give Proof.  If $V^s$ is a $j$-fold branched cover over
$\Delta(\epsilon)$, then it is homologous to $j$ times the class of
$\{0\}\times\Delta(\epsilon)$ in
$H_2(\Delta^2(\epsilon),\Delta(\epsilon)\times\partial\Delta(\epsilon))$.
Similarly, $V^u$ is homologous to $k$ times the
class of $\Delta(\epsilon)\times\{0\}$ in
$H_2(\Delta^2(\epsilon),\partial\Delta(\epsilon)\times\Delta(\epsilon))$.
Thus the intersection number of the classes $[V^s]$ and $[V^u]$ is $jk$ times
the intersection number of $\{0\}\times\Delta(\epsilon)$
and $\Delta(\epsilon)\times\{0\}$, which is 1.

For $q\in J\cap\Delta^2(\delta)$, we let $j=j_q$ denote the branching degree of
$\pi_u:V^u(q,\epsilon)\to\Delta(\epsilon)$.  Let us take a sequence $q_k\to
p$ such
that $j=j_{q_k}$ is constant and $\psi^u_{q_k}\to\psi^u\in\Psi^u_p$. Let
$\omega_k\subset\C$ denote the connected component of
$\psi^{-1}_{q_k}(V^u(q_k,\epsilon))$ containing $0$.   For each
$x_0\in\Delta(\epsilon)$ and each $k$ we have
$\#\{\zeta\in\omega_k:\pi_u\circ\psi^u_{q_k}(\zeta)=x_0\}=j$.  By [BS8,
Lemma 2.1] there
exists $r>0$ such that $\omega_k\subset\{|\zeta|<r\}$ for all $k$.  It
follows that
$\#\{|\zeta|\le r:\pi_u\circ\psi^u(\zeta)=x_0\}\ge j$. By (3.2) we have a
holomorphic map $\pi_u^{-1}:\Delta(\epsilon)\to V^u(p,\epsilon)$, so we
conclude
that $\pi_u^{-1}\pi_u\psi_p=\psi_p$ is at least $j$-to-1.  It follows from
Proposition 2.3 that $j\le2$. \qed

The sets $S:=\Delta^2(\epsilon)\cap\R^2$ and
$S_0:=\Delta^2(\delta)\cap\R^2$ are
squares in $\R^2$.  We define the {\it vertical boundary }
$\partial_v S$ (resp.\ the {\it horizontal boundary } $\partial_h S$) as the
portion of (the square) $\partial S$ which is vertical (resp.\ horizontal).
For $q\in J\cap S_0$, we define $\gamma_q^s$ as the intersection
$V^s(q,\epsilon)\cap\R^2$.  Thus $\gamma^s_q=\gamma^s_q(\cB)$ is determined
by $\cB$.  We define $\Gamma^s=\Gamma^s(\cB)$ to be the set of curves
$\gamma^s_q$
with $q\in S_0=S_0(\cB)$.  We define $\Gamma^s_j=\Gamma^s_j(\cB)$ as the set of
curves $\gamma^s\cap V^s$ with $V^s\in\cV^s_j$.   The layout of this
configuration
is illustrated in Figure 3.1:
$\gamma^s_p\in\Gamma^s_1(\cB)$, and
$\gamma^s_q,\gamma^s_r\in\Gamma^s_2(\cB)$.
By the reality condition, $\gamma^{s/u}_p\in\Gamma^{s/u}$ is a
one-dimensional set,
so $\gamma_p^{s/u}$ is regular if and only if $V^{s/u}_\epsilon(p)$ is
regular.

\epsfxsize1.8in
\centerline{\epsfbox{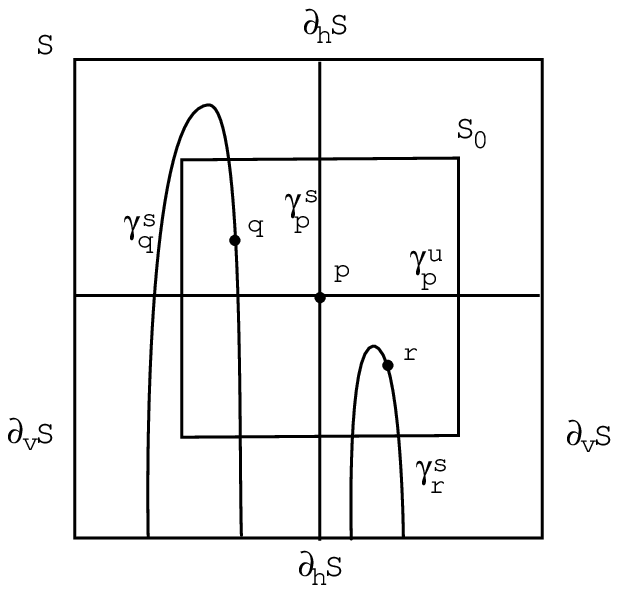}}
\centerline{Figure 3.1}
\proclaim Corollary 3.2.  If $\gamma^s\in\Gamma_j^s$ and
$\gamma^u\in\Gamma_k^u$,
then the number of points of $\gamma^s\cap\gamma^u$ counted with
multiplicity, is equal to $jk$.

\give Proof.  This is a direct consequence of Lemma 3.1 and the fact that
$V^s\cap V^u\subset\R^2$.  \qed

If $\psi\in\Psi^u_p$ has order 2, and if $\gamma^u_p$ is regular, then by
Proposition 2.6 there is an embedding $\phi$ such
that $\psi(\zeta)=\phi(\zeta^2)$.  It follows that
$\psi(\C)\cap\R^2=\psi(\R)\cup\psi(i\R)$.  Working inside a box $\cB$, we write
$(\gamma^u_p)^r:=\psi(\R)\cap\gamma^u_p$ and
$(\gamma^u_p)^i:=\psi(i\R)\cap\gamma^u_p$.  The phantom gray region
$(\gamma^u_p)^i$, as in Figure 2.1, is disjoint from $J$.  We state this
observation as follows.

\proclaim Lemma 3.3.  If $p\in J_{*,2}$ and $\gamma^u_p$ is regular, then
$p$ is
$u$ one-sided; if $p\in J_{2,*}$ and $\gamma^s_p$ is regular, then $p$ is
$s$ one-sided.

Let $q$ be a point of $J$ for which $V^s_q$ is regular For
$\psi\in\Psi^s_q$, we define the set
$\omega_\psi$ as the connected component of $\psi^{-1}V^s(q,\epsilon)$
containing
the origin.  Since $\pi_u\circ\psi$ is an entire function, $\omega_f$ is simply
connected.  By the reality condition on $\psi$, $\omega_\psi$ is invariant
under
complex conjugation.  Thus $\omega_\psi\cap\R$ is a (connected) interval
$(-a,b)$.  It follows that $\gamma_q^s$ is connected submanifold of $S$.
We refer
to $\psi(-a)$ and $\psi(b)$ as the endpoints of $\gamma^s_1$.  Since
$\partial
V^s(q,\epsilon)\subset\Delta(\epsilon)\times\partial\Delta(\epsilon)$, it
follows that the endpoints of $\gamma^s_q$ lie in $\partial_hS$.

\proclaim Lemma 3.4.  If
$\gamma\in\Gamma^s$, then $\gamma\in\Gamma^s_1$ if the
endpoints of $\gamma$ lie in different components of $\partial_hS$.
Otherwise (if the endpoints lie in the same component of $\partial_hS$),
$\gamma\in\Gamma^s_2$.

\give Proof.  The horizontal boundary $\partial_v
S$ consists of fibers of the projection $\pi_s$ intersected with $\R^2$.
By Lemma 3.1, the multiplicity of the projection is no greater than 2. If
$\gamma^s$ intersects one of the fibers in two points, then the multiplicity
is in fact equal to 2.

Now we prove the first assertion of the Lemma.  Suppose that $\gamma^s$ has
one endpoint in each component of $\partial_hS$.  Then for each point
$t\in\Delta(\epsilon)\cap\R$, there is a point $s\in(-a,b)$ such that
$\psi(s)=t$.  Now there cannot be a point $\zeta\in\C-\R$ with
$\psi(\zeta)=t$, for by the reality condition we would have
$\psi(\bar\zeta)=t$, which would give 3 solutions.  Finally, we cannot have the
situation where $\pi^{-1}_s(t)\cap\gamma^s$ consists of exactly two points.
For,
in this case, we may assume that $\pi_s(\psi(-a))=-\epsilon$ and
$\pi_s(\psi(b))=\epsilon$.  Then, arguing as in Calculus, there must be a
nearby
$t'\in(-\epsilon,\epsilon)$ for which $\pi^{-1}_s(t')$ consists of 3 points.
\qed

\proclaim Lemma 3.5.  Let $\cB$ be a regular box.  After possibly shrinking
$\delta$, it follows that if $q\in J_{2,*}\cap S_0$ and if $\gamma^s_q$ is
regular,
it follows that $\gamma^s_q\in\Gamma^s_1$.  Similarly, if $q\in J_{*,2}\cap
S_0$,
and if $\gamma^u_q$ is regular, then $\gamma^u_q\in\Gamma^u_1$.

\give Proof.  It follows from Lemma 3.1 and Corollary 3.2 that $\gamma^s_q$
belongs
to $\Gamma^s_1$ or $\Gamma^s_2$.  If there is no $\delta$ satisfying the
conclusion of the Lemma, then there is a sequence $q_j\to0$, $q_j\in J_{2,*}$,
$\gamma^s_{q_j}$ regular, and $\gamma^s_{q_j}\in\Gamma^s_2$.  Since
$\gamma^s_{q_j}$
is regular, there exists $\psi^s_j\in\Psi^s_{q_j}$ and a holomorphic embedding
$\phi_j$ such that $\psi^s_j(\zeta)=\phi_j(\zeta^2)$.  We may extract a
subsequence
such that there is a limit $\psi^s_j\to\psi\in\Psi^s_p$.

Now since $\psi^s_{q_j}\in\Gamma^s_2$ it follows that for
$x_0\in\Delta(\epsilon)$,
$\pi_s^{-1}(x_0)\cap V^s_{p_j}$ consists of two points.  Thus
$$\#(\pi_s\circ\psi^s_j)^{-1}(x_0)=\#\{\zeta\in\C:\pi_s\phi(\zeta^2)=x_0\}\ge4.
$$
By [BS8, Lemma 2.1] this set is contained in a disk $\{|\zeta|<r\}$,
independent of
$j$.  Letting $j\to\infty$ we obtain $\#(\pi_x\circ\psi)^{-1}(x_0)\ge4$.   By
(3.2), there exists $a_0\in V^s_p$ such that
$V^s_p\cap\pi_s^{-1}(x_0)=\{a_0\}$ is a single point so
$\#(\pi_s\psi)^{-1}(x_0)=\#\psi^{-1}(a_0)\ge4$, which contradicts Lemma 2.3.
\qed
\epsfxsize4in
\centerline{\epsfbox{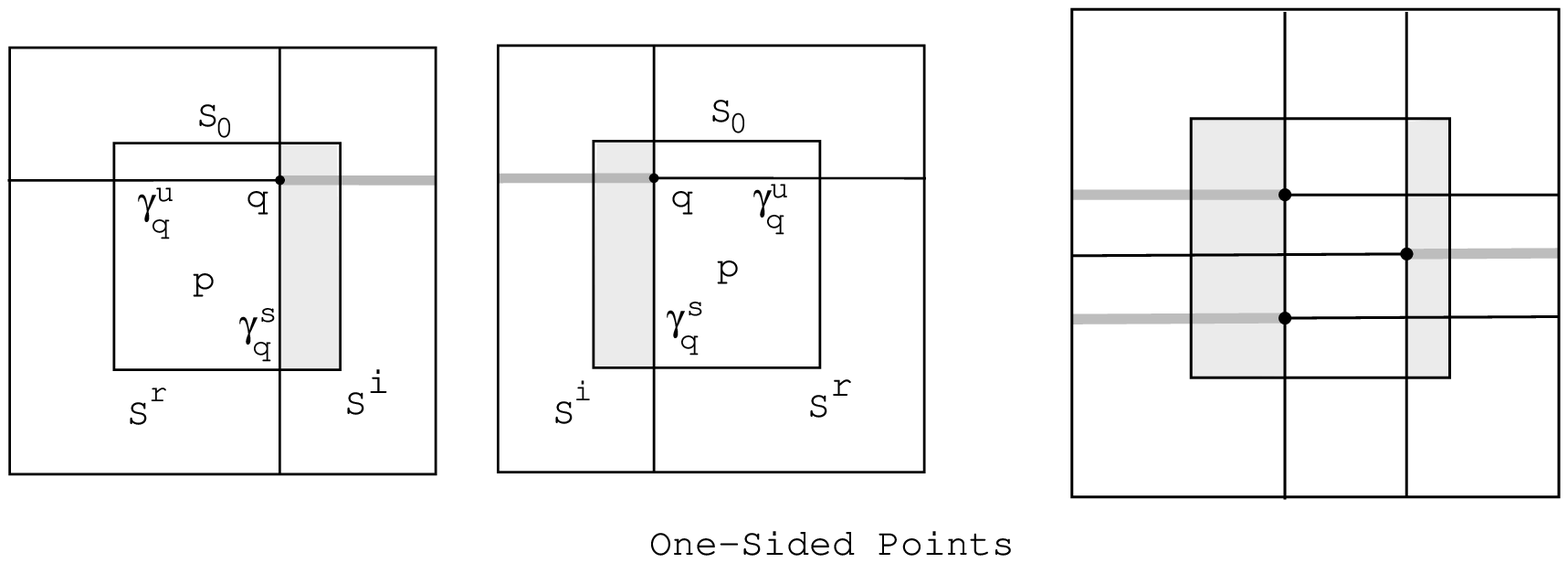}}
\centerline{Figure 3.2}
Let $\cB$ be  a regular box, and let $q\in S_0\cap J_{*,2}$ be a point with
$\gamma^s_q\in\Gamma^s_1$.  Then $S-\gamma^s_q$ consists of two components,
which
we may label $S^r$ and $S^i$, as in Figure 3.2.  That is, $S^r$ contains
the variety
of $\psi(\R)$ at $q$, and $S^i$ contains the local variety of the phantom
region
$\psi(i\R)$ at $q$.

\proclaim Lemma 3.6.  For $p\in J_{2,2}$, we may construct a regular box $\cB$
about $p$.  For $q\in S_0\cap J_{2,2}$, $\gamma^s_q$ belongs to
$\Gamma^s_1$. If
we split $S-\gamma^s_q=S^r\cup S^i$ as above, then $S_0\cap S^i\cap
J=\emptyset$.  (And the corresponding statement holds for $S-\gamma^u_q$.)

\give Proof.  By hypothesis, $J_{2,2}\ne\emptyset$.  Thus
$2=\max\tau^u=\max\tau^s$, so by [BS8, Proposition 5.2] $J_{2,2}$ is a
hyperbolic
set.  It follows that $\gamma^{s/u}_q$ are nonsingular and transversal.  In
particular, for $q=p$, we may construct a regular box $\cB$ about $p$.  If
$\cB$
is sufficiently small, then it follows by hyperbolicity that
$\gamma^{s/u}\in\Gamma^{s/u}_1$ for all $q\in S_0\cap J_{2,2}$.

To complete the proof, we must show that $S_0\cap S^i\cap
J=\emptyset$.  For otherwise, if there exists $r\in S_0\cap S^i\cap J$, then by
Corollary 3.2 $\gamma^s_r\cap\gamma^u_q\ne\emptyset$.  Since $\gamma^s_r$
cannot
intersect $\gamma^s_q$, it follows that $\gamma^s_r\cap\gamma^u_q$ must lie
inside the phantom region of $\gamma^u_q$, which is forbidden.\qed

\proclaim Theorem 3.7.  $J_{2,2}$ is finite.

\give Proof.  By Lemma 3.6 we may construct a regular box $\cB$ about any $p\in
J_{2,2}$.  Let us select a finite family of boxes $\cB$ such that the
corresponding sets $S_0$ cover $J_{2,2}$.  Let us fix one of these sets $S_0$.
If $q\in J_{2,2}\cap S_0$, then $q$ corresponds to one of the four
types of doubly one-sided point pictured on the left hand side of Figure 3.3.
For each of these four cases, it follows from Lemma 3.6 that the set $J\cap S$
must lie in the quadrant bounded by the solid lines.

\epsfxsize4in
\centerline{\epsfbox{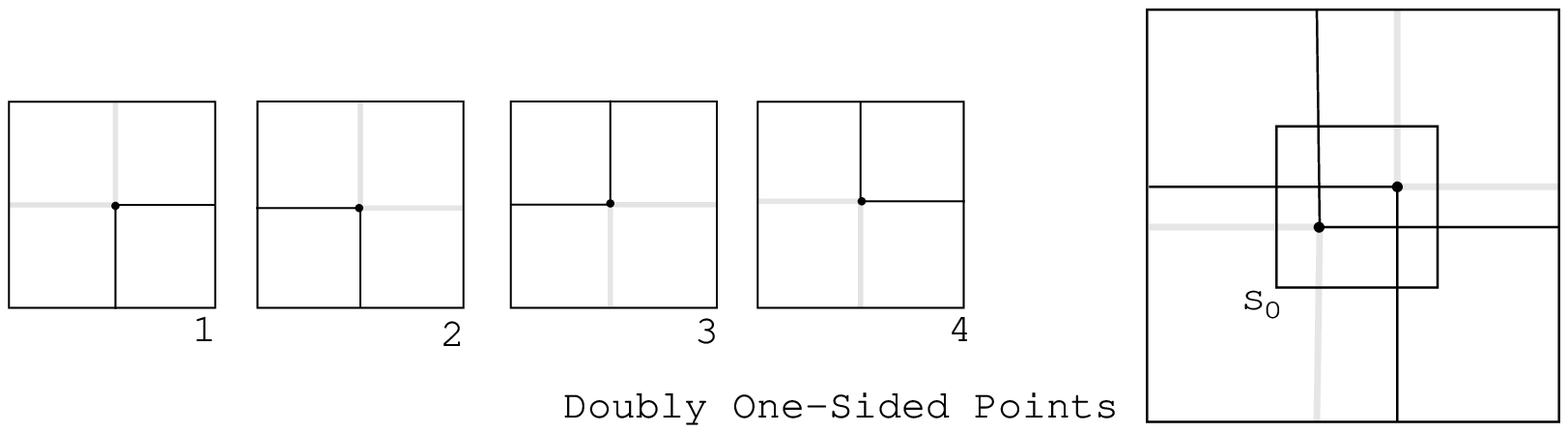}}
\centerline{Figure 3.3}

Now we consider the possibility that $S_0\cap J_{2,2}$ might consist of
more than
one point.  Let us start by supposing that a pair of points $p_1,p_2\in J$
belong to
the same box $S_0$.  The only way that two types of box can both occupy the
same
set
$S$ is if they are of type 1 and 3 or type 2 and 4.  The situation where
points of
type 2 and 4 occupy the same box $S$ is pictured on the right hand side of
Figure
3.3, and by Lemma 3.6 $J\cap S_0$ is contained in the shaded region.  It
follows
that $J\cap S_0$ cannot contain a third point $r$.  If there were, it would
lie in
the shaded portion; but as $r$ is one-sided, then the phantom (gray prong)
region
would necessarily intersect the sides of the shaded region.

Since we have covered $J_{2,2}$ by finitely many sets $S_0$ and
$\#(J_{2,2}\cap S)\le2$, it follows that $J_{2,2}$ is finite.  \qed

\proclaim Theorem 3.8.  If $p\in J_{2,1}\cup J_{1,2}$, then $V_p^{s/u}$ is
regular
at $p$.

\give Proof.  Without loss of generality we may assume that $p\in J_{1,2}$.  If
$\alpha(p)\cap J_{1,2}\ne\emptyset$ it follows from [BS8, Theorem 5.5] that
$V^u_\epsilon(p)$ is regular.  The other possibility is that $\alpha(p)\subset
J_{2,2}$ by Corollary 2.4.  Since $J_{2,2}$ is hyperbolic and finite, it is an
isolated hyperbolic set. By [R, p.\ 380] there exists $q\in J_{2,2}$ such that
$p\in W^u(q)$.  By (1.4) it follows that $V^u_p\subset W^u(q)$, and so
$V^u_\epsilon(p)$ is regular.  \qed

For $p\in J_{1,2}$ it follows from Theorem 3.8 that $\gamma^u_p$ is
regular.  By
Theorem 2.7, $\gamma^s_p$ is regular and transverse to
$\gamma^u_p$.  If $p\in J_{2,2}$, then we saw in the proof of Theorem 3.7 that
$\gamma^s_p$ and $\gamma^u_p$ are regular and transverse at $p$.  Thus if $p\in
J_{*,2}$ we may construct a regular box
$\cB$ centered at $p$.  The following result will involve shrinking this box
$\cB$.  Before giving the proof, we make an observation concerning the
relationship between shrinking and the multiplicities of varieties.    Let
$\pi_s$ and $\pi_u$ be the projections associated with $\cB$.  If
$V\in\cV^u_m$, then the projection
$\pi_u:V\to\Delta(\epsilon)$ has mapping degree $m$.  This is equivalent to the
statement that the total multiplicities of the critical points of $\pi_u$ is
$m-1$.  It follows that if we shrink the box $\cB$ to
$\cB':=\{q\in\cB:\pi_u(q)\in\Delta(\epsilon'_1),
\pi_s(q)\in\Delta(\epsilon'_2)\}$ for some $\epsilon'<\epsilon$, then each
component of $V\cap\cB'$ belongs to $\cV^u_j(\cB')$ for $j\le m$.  It follows
that if $S'=\R^2\cap\cB'$ is a regular box obtained by shrinking $\cB$ in this
way, then for each $\gamma\in\Gamma^u_1(S)$, $\gamma\cap S'\in\Gamma^u_1(S')$.
And for each $\gamma\in\Gamma^u_2(S)$, we have that either $\gamma\cap S'$ is
connected and belongs to $\Gamma^u_2$; or $\gamma\cap S'$ consists of two
components, each of which belongs to $\Gamma^u_1(S')$.  In this sense,
$\Gamma^{s/u}_1(S)$ is preserved under shrinking.

\proclaim Lemma 3.9.  For $p\in J_{*,2}$ there is a regular box $\cB$
centered at
$p$ with the properties:
\item{1.} For all $\gamma\in\Gamma^u_2$, $\bar\gamma\cap\partial_v S\subset\bar
S^r$.
\item{2.}  For all $\gamma\in\Gamma^u_2$, $\gamma\cap S^r$ consists of two
components $\gamma_1$ and $\gamma_2$, as in Figure 3.4.
\item{3.} For all $\eta\in\Gamma^s_2$,
$\#(\gamma_1\cap\eta)=\#(\gamma_2\cap\eta)=2$.

\give Proof.  As we noted in the previous paragraph, we may construct a regular
box  $\cB$ centered at $p$.  Now we show that we may shrink $S_0$ and $S$
so that
1, 2, and 3 hold.  Since $p\in J_{*,2}$, $\gamma^u_p$ will have a phantom
region, which we may assume extends to the right hand side, as in Figure 3.4.
In order to establish 1, we must show that for $\gamma\in\Gamma^u_2$, the
endpoints of $\gamma$ lie in the left hand side of the vertical boundary of
$S$. Let $\gamma'\in\Gamma^u_2$ denote any unstable arc whose endpoints lie in
the right hand side of $\partial_vS$.  For $r\in S_0\cap J$,
$\gamma^s_r\cap\gamma$ consists of two points, which means that any $\gamma'$
must loop around to the left of $\gamma^s_r$.  If we shrink $S$ in the unstable
direction, i.e., replace it with $S\cap\pi^{-1}\Delta(\epsilon')$ with
$\epsilon'>0$ small enough that there exists $r\in S_0\cap J$ with
$\gamma^s_r\cap S'=\emptyset$, then all $\gamma'$ become simple in $S'$. That
is, $\gamma'\cap S'\in\Gamma^u_1(S')$.

Assertion 2 follows from assertion 1, as is illustrated in the left hand
side of
Figure 3.4.

\epsfxsize3.8in
\centerline{\epsfbox{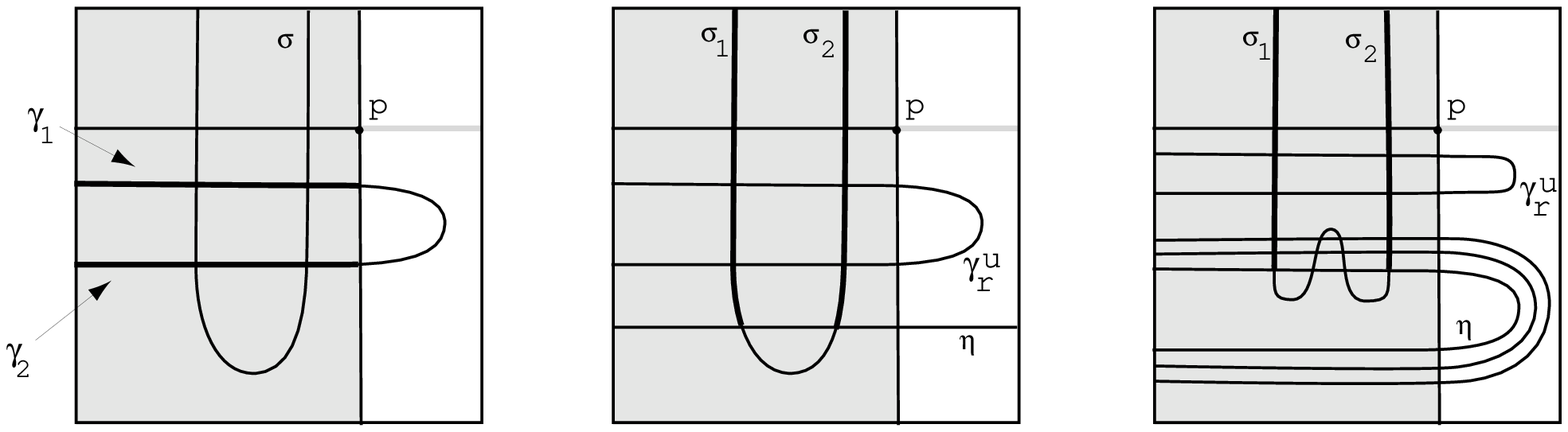}}
\centerline{Figure 3.4}

To prove assertion 3, we consider first the case where there is an
$\eta\in\Gamma^u_1$ lying below $\gamma^u_p$, as in the central picture in
Figure
3.4.  (The case where there
$\eta$ lies above $\gamma^u_p$ is analogous.)  We may shrink
$S_0$ so that for all $r\in J\cap S_0$, $\gamma^u_r$ lies between $\eta$ and
$\gamma^u_p$.  In this case we consider
$\gamma\in\Gamma^u_2$ lying between
$\eta$ and
$\gamma^u_p$ and $\sigma\in\Gamma^s_2$.  The case drawn in the central
picture in
Figure 3.4 shows the endpoints of $\sigma$ in the top of the boundary
$\partial_hS$.  (The other case, where the endpoints are in the bottom
portion of
$\partial_hS$ is analogous.)  As is pictured, $\eta$ cuts off two pieces
$\sigma_1$ and $\sigma_2$, and each of these intersects $\gamma_j$, $j=1,2$.
Thus $\#(\sigma\cap\gamma_j)=2$.  This proves assertion 3 in this case.

The alternative to this case is that $\eta\in\Gamma^u_2$ for all curves
$\eta\in\Gamma^u$ lying below $\gamma^u_p$.  If this happens, we consider
$G(\eta)$, which is the set of all $\gamma^u_r$ such that
$\gamma^u_r\in\Gamma^u_2$, and $\gamma^u_r$ separates $\eta$ from $p$.  We
claim
that we may shrink $S_0$ such that for $r\in S_0$ below $\gamma^u_p$,
$\gamma^u_r$ lies between $G(\eta)$ and $\gamma^u_p$.  If this happens, then we
see that any $\sigma\in\Gamma^s_2$ has two components $\sigma_1$ and $\sigma_2$
as in the right hand side of Figure 3.4.  These components intersect
$\gamma$ as
desired.  The alternative is that there are points $r_j\in S_0\cap J$,
lying below $\gamma^u_p$, and such that $r_j\to p$.  But then we have that
$\gamma^u_{r_j}\to\gamma^u_p$ in the topology of the Hausdorff metric.  In this
case, we let $S'$ denote an arbitrarily small shrinking of $S$, and it follows
that $\gamma^u_{r_j}\cap S'$ is ultimately disconnected.  Thus a component of
$\gamma^u_{r_j}$ serves as the curve $\eta$ as we considered at first.
\qed

\proclaim Theorem 3.10.  $J_{2,1}\cup J_{1,2}$ is finite.

\give Proof.  It suffices to show that $J_{1,2}$ is finite.  Write
$X=J_{1,2}\cup J_{2,2}$.  For each
$p\in X$, we may
construct a regular box $\cB$ centered at $p$, satisfying the conclusions of
Lemma 3.9.   Since $X$ is compact, we may select a finite number of regular
boxes $\cB$ such that the sets $S_0$ cover $X$.

Let us fix one of these boxes.
We claim: {\sl There are (at most) two verticals, $\gamma^s_{p}$ and
$\gamma^s_{q}$, with the property that
$S_0\cap X\subset\gamma^s_{p}\cup\gamma^s_{q}$.}  Since $p\in
X$, it is $u$ one-sided.  Without loss of generality we may assume
that the phantom region of $\gamma^u_p$ is on the right, so that
$S_0\cap J$ lies to the left of $\gamma^s_p$.  To establish the claim, we show
that for any two points
$q,r\in X\cap S_0$ such that $q,r\notin \gamma^u_p$, then it follows that
$\gamma^s_q=\gamma^s_r$.

Let us assume first that $\gamma^u_r$ and $\gamma^u_q$ both belong to
$\Gamma^u_1$. If $r\notin\gamma^s_q$, then we have
$\gamma^s_q\cap\gamma^u_r\ne\emptyset$, as pictured in Figure 3.5.  There are
three possibilities for the $\gamma^s_r$.  The first (on the left of Figure
3.5)
is that $\gamma^s_r\in\Gamma^s_1$.  But this is not possible, since the phantom
region of $\gamma^u_q$ blocks $\gamma^s_r$ from reaching the upper portion of
$\partial_hS$.  The next two possibilities are that $\gamma^s_r\in\Gamma^s_2$.
In both cases, $\gamma^s_r\cap\gamma^s_q=\emptyset$, and $\gamma^s_r$ cannot
intersect the phantom region of $\gamma^u_q$.  In the central picture of Figure
3.4, we see that $\gamma^s_r$ must come out of the box bounded by $\gamma^u_q$
and $\gamma^s_q$.  But the portion that is drawn shows
$\#(\gamma^s_r\cap\gamma^u_r)=2$, and thus $\gamma^s_r$ cannot intersect
$\gamma^u_r$ again.  Thus $\gamma^s_r\cap\gamma^u_q=\emptyset$, which is a
contradiction.  In the last case, on the right of Figure 3.5, we again have
$\#(\gamma^s_r\cap\gamma^u_r)=2$ already, so it is not possible for
$\gamma^s_r$
to cross $\gamma^u_r$ again, so it cannot intersect $\partial_hS$, which is a
contradiction.

Now let us suppose that one or both of $\gamma^u_r$, $\gamma^u_q$ belongs to
$\Gamma^u_2$.  We will refer to these as $\gamma^u_x$.  As in Lemma 3.9,
$\gamma^u\cap S^r$ consists of two pieces, $\gamma_1$ and $\gamma_2$.  One of
these, say  $\gamma_1$, contains the point $x$.  By Lemma 3.9, we have
$\#(\gamma^s\cap\gamma^u_x)=2$ for any $\gamma^s\in\Gamma^s_2$.  Thus we
replace
$\gamma^u_x$ in Figure 3.5 by $\gamma_1$ and proceed as before.  This
proves the
claim.

Since there are only finitely many regular boxes in our covering of $X$, it
follows that $J_{1,2}$ is contained in the union of a finite set
$\{\gamma^s_1,\dots,\gamma^s_N\}$ of segments of stable manifolds.  We let
$\Gamma^s_j$ denote the closure of the union of all the arcs
$\gamma^s_i$ which are contained in the global stable manifold
$W^s(r_j)\supset\gamma^s_j$.  It follows that $f$ permutes the finite family of
sets
$\{\Gamma^s_1,\dots,\Gamma^s_N\}$.  Thus, passing to a power of $f$, we have
$f^t(\Gamma_1^s\cap J_{1,2})=\Gamma^s_1\cap J_{1,2}$.  Now $\Gamma^s_1\cap
J_{1,2}$
is an $f$-invariant, compact subset of $W^s(r_1)$, and $f^t$ is contracting on
$W^s(r_1)$, so
$\Gamma^s_1\cap J_{1,2}$ must be a single point.  We conclude that $J_{1,2}$ is
finite. \qed
\epsfxsize3.8in
\centerline{\epsfbox{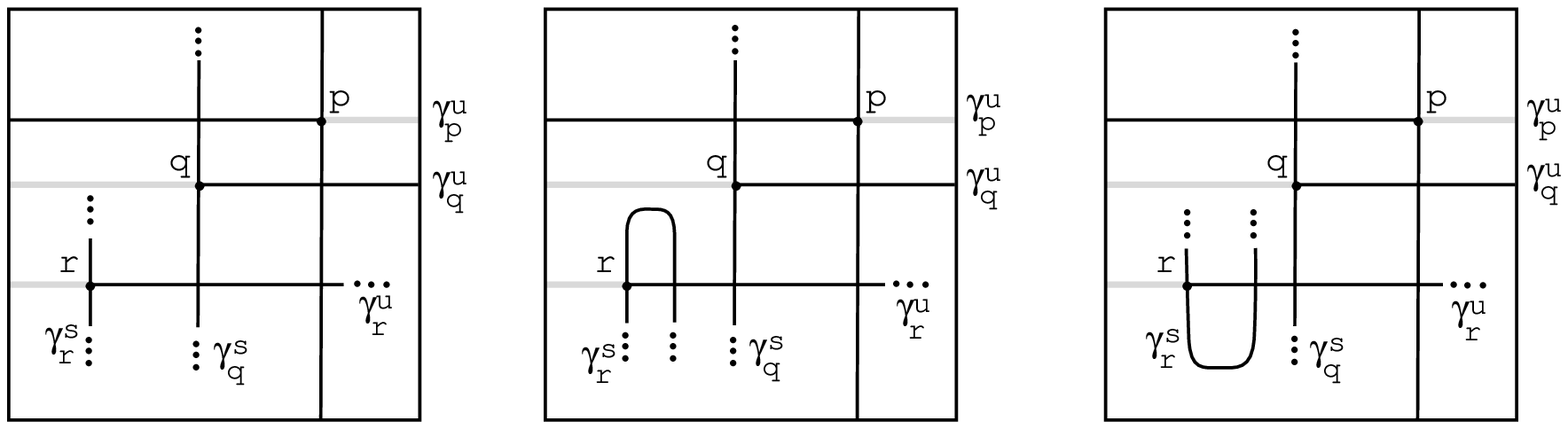}}
\centerline{Figure 3.5}

\proclaim Corollary 3.11.  For any $p\in J$ and $\psi\in\Psi_p$, $\psi(\C)$
is a
nonsingular (complex) submanifold of $\C^2$, and $\psi(\C)\cap\R^2$ is a
nonsingular (real) submanifold of $\R^2$.

\give Proof.  If $\psi$ has no critical point, then $\psi(\C)$ is nonsingular.
And by our earlier discussion of the reality condition, it follows that if
$\psi$ has no critical point, then $\psi(\C)\cap\R^2$ is a nonsingular, real
one-dimensional submanifold of $\R^2$.

If $\psi\in\Psi^u_p$ has a critical point, then by Proposition 2.5, $\psi$
there
is just one critical point $\zeta_0$.  The sequence $\tilde
f^{-n}\psi=f^{-n}\psi(\lambda(p,-n)^{-1}\zeta)\in\Psi_{f^{-n}p}$ has a critical
point at $\lambda(p,{-n})\zeta_0$.  For a subsequence
$f^{-n_j}p\to q\in\alpha(p)$, we may pass to a further subsequence such that
$\tilde f^{-n_j}\psi$ converges to $\hat\psi\in\Psi_q$.  Since
$\lambda(p,-n)\to0$ as $n\to\infty$,  $\hat\psi$ has a critical point at
$\zeta=0$.  Thus $q\in J_{*,2}$.  Since $J_{2,*}\cup J_{*,2}$ is a finite
set of
saddle points, we have $p\in W^u(p)$.  Thus $\psi(\C)\subset W^u(q)$, and
so this
set and $\psi(\C)\cap\R^2$ are both regular.  \qed

\section 4.  Hyperbolicity and Tangencies

In \S3 we showed that $\cC:=J_{2,*}\cup J_{*,2}$ is a finite union of saddle
points.  We show next that all tangential intersections lie in stable manifolds
of $J_{*,2}$ and unstable manifolds of $J_{2,*}$.  In Theorem 4.2 we show that
for $p\in J_{*,2}$, the stable manifold $W^s(p)$ contains a heteroclinic
tangency.  The condition for hyperbolicity is characterized (Theorem 4.4) in
terms of the existence of heteroclinic tangencies.

\proclaim Theorem 4.1.  If $p\in J$ and $E^s_p=E^u_p$, then $p\in
W^s(J_{*,2})\cap W^u(J_{2,*})$.

\give Proof.  If $p$ is a point of tangency, then by Theorem 2.7,
$\alpha(p)\subset J_{2,*}$ and $\omega(p)\subset J_{*,2}$.  Since $J_{2,*}\cup
J_{*,2}$ is a finite set of saddle points, it follows that there exist $q\in
J_{2,*}$ and $r\in J_{*,2}$ such that $p\in W^u(q)\cap W^s(r)$. \qed

\proclaim Theorem 4.2.  If $p\in J_{*,2}$ then there exists $q\in J_{2,*}$
such that $W^s(p)$ intersects $W^u(q)$ tangentially.

\give Proof.  By \S3, $\cC$ is finite. For each point $p\in J_{*,2}$, there
is a
point $q\in\cC$ such that $W^s(p)$ intersects $W^u(q)$ tangentially, by [BS8,
Theorem 8.10].  And by Theorem 2.7 we have $q\in J_{2,*}$.  \qed

\proclaim Corollary 4.3.  If $J_{1,2}\ne\emptyset$, then $J_{2,1}\ne\emptyset$.

\proclaim Theorem 4.4.  The following are equivalent for a real, polynomial
mapping of maximal entropy:
\item{1.}  $f$ is not hyperbolic.
\item{2.}  $J_{2,*}\cup J_{*,2}$ is non-empty.
\item{3.}  There are saddle points $p$ and $q$ such that $W^s(p)$ intersects
$W^u(q)$ tangentially.

\give Remark.  By Theorem 4.1, the saddle points $p$ and $q$ in condition 3
satisfy $p\in J_{*,2}$ and $q\in J_{2,*}$.

\give Proof.  $(1)\Rightarrow(2)$.  If $J_{2,*}\cup J_{*,2}\ne\emptyset$,
then by
Theorem 4.2 there is a tangency between $W^s(J_{*,2})$ and $W^u(J_{2,*})$. Thus
$f$ is not hyperbolic.

$(2)\Rightarrow(1)$.  If $J_{2,*}\cup J_{*,2}=\emptyset$, then $J=J_{1,1}$.  It
follows that $J_{1,1}$ is compact, so by [BS8, Proposition 5.2] $J_{1,1}$
is a hyperbolic set.

The implication $(2)\Rightarrow(3)$ follows from Theorem 4.2, and
$(3)\Rightarrow(2)$ follows from Theorem 2.7.\qed

Let $\cT$ denote the set of points of tangential intersection between $W^s(a)$
and $W^u(b)$, for $a,b\in J$.  By [BS8, Theorem 8.10], $\cT$ is a discrete
subset of $J_{1,1}$.  Since the parametrizations are nonsingular in $J_{1,1}$,
the curves $\cW^{s/u}=\{\gamma^{s/u}_r:r\in J_{1,1}\}$ form a lamination of a
neighborhood of $p$.  If $p\in J_{1,1}-\cT$, the the laminations
$\cW^s$ and $\cW^u$ are transverse at $p$, so they define a local product
structure on $J$ in a neighborhood of $p$.

Let us fix a point $p\in\cT$.  We cannot construct a regular box centered
at $p$
since it is a point of tangency, but we will construct a box with many of the
same properties.  We choose a real analytic coordinate system such that the
square $S:=\{|x|,|y|<\epsilon\}\subset\R^2$ is centered at $p=(0,0)$, and
has the
properties that $\gamma^s_p=\{x=0,|y|<\epsilon\}$, and the projections
$\pi_u:\gamma^u_p\to(-\epsilon,\epsilon)$ are proper, where as before
$\pi_u(x,y)=x$, $\pi_s(x,y)=y$, and we use the notation
$\gamma^{s/u}_q:V^{s/u}_q\cap S\cap\R^2$.  By Theorem 2.8, the multiplicity of
the intersection of $\gamma^s_p$ and $\gamma^u_p$ at $p$ must be 2.  Thus
$\gamma^u_p\in\Gamma^u_2(S)$, and so by Lemma 3.4 $\gamma^u_p$ lies to one side
of $\gamma^s_p$, as is pictured on the left hand side of Figure 4.1.  For
$S_0:=\{|x|,|y|<\delta\}$ sufficiently small, we have (3.1), and
$\pi_{s/u}:\gamma^{s/u}_q\to(-\epsilon,\epsilon)$ is proper for $q\in
S_0\cap J$.

The configuration of the curves in the third picture of Figure 4.1 follows from
Lemma 3.4 and Corollary 3.2, since there must be two points of intersection
between stable and unstable manifolds in $S$.  This arrangement is associated
with the failure of topological expansivity.
\medskip
\epsfxsize4.5in
\centerline{\epsfbox{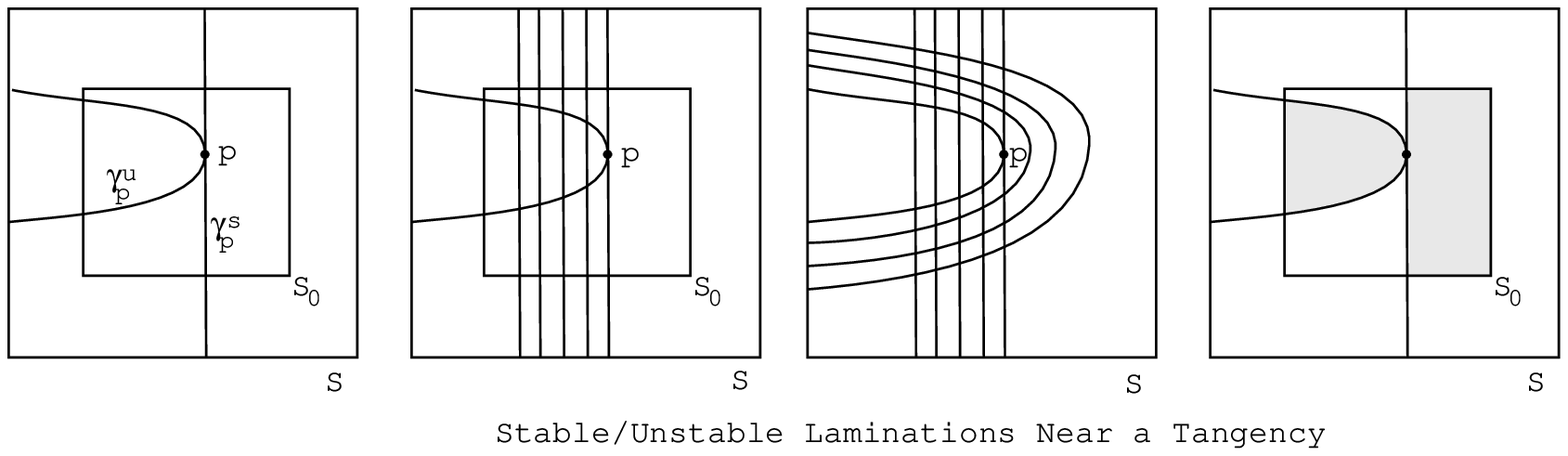}}
\centerline{Figure 4.1}

\proclaim Corollary 4.5.  If $r\in\cT$, then there is a neighborhood $S_0$
of $r$ such that $J\cap S_0$ is disjoint from the region shaded in the
right hand
picture in Figure 4.1.

\section 5.  One-Sided Points

We have shown that the set of critical points $\cC$ is a finite set of
one-sided
points.  We use one-sided points to show (Theorem
5.2) that $K$ is always a Cantor set.  We analyze more carefully  the
possibilities for one-sided points, and we obtain Propositions 5.8 and 5.9,
which combine to prove Theorem 3 in the Introduction.\footnote*{We
wish to thank Andr\'e de Carvalho for a suggestion that resulted in this
part of the paper.}

\proclaim Theorem 5.1.  If  $f$ is hyperbolic, then there exist stably and
unstably one-sided points.

\give Proof.  The set $K$ is saturated in the sense that $W^s(p)\cap
W^u(q)\subset K$ for all $p,q\in K$.  We use the following result of
Newhouse and
Palis [NP], as it is presented in [BL, Proposition 2.1.1, item 6]: If $f$ does
not have an unstably one-sided point, then $K$ is a hyperbolic attractor. Since
$K$ is a basic set, it follows that if it is an attractor, then the set of
points attracted to $K$ is open.  On the other hand, the set of points
attracted
to $K$ is $K^+$, which is a closed, proper subset of $\R^2$, and is thus not
open.  Thus $f$ has an unstably one-sided point. Repeating the argument for
$f^{-1}$ gives a stably one-sided point. \qed

\proclaim Theorem 5.2.  If $f$ has maximal entropy, then $K$ is a Cantor set.

\give Proof.  Since $K$ is the zero set of a continuous, plurisubharmonic
function $G$ on $\C^2\supset\R^2$, it follows that no point of $K$ can be
isolated.  Thus it suffices to show that $K$ is totally disconnected.  Both
$W^s(K)$ and $W^u(K)$ are laminations in a neighborhood of
$\cJ_{1,1}=K-\cC$.  Let $\cT$ denote the tangencies between $W^s(K)$ and
$W^u(K)$.  By [BS8, Theorem 8.10] $\cT\cup\cC$ is a countable, closed set.
Thus it suffices to show that $K$ is totally disconnected in a neighborhood of
each point of $K-(\cT\cup\cC)$.

Now each point of $K-(\cT\cup\cC)$ has a neighborhood $R$ such that $R\cap K$
has local product structure.  The local product structure means that for any
$r\in K\cap U$, $R\cap K$ is homeomorphic to $(K\cap W^u_{R,loc}(r))\times
(K\cap W^s_{R,loc}(r))$.

By Theorem 5.1, there are an $s$ one-sided periodic point $p$ and a
$u$ one-sided
point $q$.  Let $A$ denote the set of transverse intersections of $W^u(p)$ and
$W^s(q)$.  By [BLS, Theorem 9.6] $A$ is dense in $K\cap U$.  It follows
from the local product structure that the transverse intersections between
$W^u(p)$ and $W^s_{R,loc}(r)$ are dense in $W^s_{R,loc}(r)$ for each $r\in
K\cap
U$.  Since $p$ is $s$ one-sided, $K\cap W^s(p)$ lies to one side of $p$ in
$W^s(p)$.  This one-sidedness propagates along the unstable manifold $W^u(p)$,
and so for any point $b\in W^u(p)\cap W^s_{R,loc}(r)$, $K\cap  W^s_{R,loc}(r)$
lies (locally) to one side of $b$ in $W^s_{R,loc}(r)$.  Thus the set of
disconnections of $K\cap  W^s_{R,loc}(r)$ is dense in $K\cap W^s_{R,loc}(r)$,
which means that $K\cap  W^s_{R,loc}(r)$ is totally disconnected.
Similarly, $K\cap W^u_{R,loc}(r)$ is totally disconnected.  By the local
product structure, $K\cap R$ is totally disconnected, and thus $K$ is
disconnected.  \qed

By a (topological) attractor we will mean a compact, invariant $S$ whose stable
set $W^s(S):=\{q:\lim_{n\to+\infty}d(f^nq,S)=0\}$ has nonempty interior.

\proclaim Corollary 5.3.  If $f$ has maximal entropy, then $K^+$ and $K^-$
have no
interior, and thus $K$ contains no attractors or repellors.

\give Proof.  As in the proof of Theorem 5.2, we consider a local
product neighborhood $R$ of a point of $K-(\cT\cup\cC)$.  By the local product
structure,  $K^+\cap R$ is homeomorphic to $(K\cap W^u_{R,loc}(r))\times
W^s_{R,loc}(r)$.  Since $K\cap W^u_{R,loc}(r)$ is totally
disconnected, it contains no interior. Thus $R\cap K^+$ contains
no interior.

If $S$ is an attractor, it must be contained in $K$, and thus the basin
$B(S)$ must be contained in $K^+$.  However, since $K^+$ has no interior, it
follows that $B(S)$ can have no interior.  \qed

Recall that if $p$ is $u$ one-sided, then $W^u(p)-\{p\}$ has a component which
is disjoint from $J$.   No point of $J\cap W^u(p)$ can be isolated, so only one
of the components of $W^u(p)-\{p\}$ can be disjoint from $J$.  We call this
component the (unstable) separatrix associated with $p$.

Now let us note that if $p$ is a saddle point of $f$ with period $n$, then
$Df^n(p)$ has eigenvalues $|\lambda^u|>1>|\lambda^s|>0$.   If $p$ is
$u$ one-sided, $f$ must preserve the unstable separatrix, and so
$\lambda^u>0$.  Similarly, if $p$ is $s$ one-sided, we have $\lambda^s>0$.

\proclaim Lemma 5.4.  Let $p\in K$ be $u$ one-sided, and let $S$ be the
separatrix associated with $p$ and which is disjoint from $K$.  Then  $S$ is
properly embedded in $\R^2-\{p\}$.

\give Proof.  Consider the uniformization $\psi^u:\C\to W^u(p)$ of the complex
unstable manifold through $p$. Since $p$ is $u$ one-sided, we may assume that
its separatrix $S$ corresponds to the positive real axis in $\C$, and thus
$G^+\psi^u(\zeta)>0$ for $\zeta\in\R$, $\zeta>0$.  Now $\psi^u(0)=p$, and there
is $\lambda^u>1$ such that $G^+\psi^u(\lambda^u)=d\cdot G^+\psi^u(\zeta)$.
Thus $\lim_{\zeta\to+\infty}G^+\psi^u(\zeta)=+\infty$.  Since $\{G^+\le c\}\cap
J^-$ is compact, it follows that $\lim_{\zeta\to+\infty}\psi^u(\zeta)=+\infty$,
and thus $S$ is properly embedded in $\R^2$.  \qed

Let $\cO=\{p_1,\dots,p_j\}$ denote the set of $u$ one-sided points.  Let
$S_1$, \dots $S_j$ denote the corresponding separatrices.  By Lemma 5.4, $S_i$
is an arc in the sphere $S^2=\R^2\cup\{\infty\}$ which connects $p_i$ to
$\infty$.  We let $\hat S_i$ denote the germ at infinity of $S_i$, i.e.\ $\hat
S_i$ denotes the set of sub-arcs of $S_i$ containing $\infty$.  We consider a
system of neighborhoods $W$ of $\infty$ in $S^2$ with the property that
$\partial
W$ is homeomorphic to $S^1$, and $\partial W\cap S_i$ consists of a unique
point,
for each $i$.  It follows that for such a neighborhood, $W\cap
(\R^2-\bigcup_{i=1}^j S_i)$ consists of $j$ open sets $V$, each containing
$\infty$ in its closure.  Let $\hat V$ denote the germ at infinity
corresponding
to $V$.  We will define a graph $\cG$ whose vertices be the classes $\hat S_i$,
and whose edges are the germs $\hat V$. A vertex $\hat S_i$ is contained in an
edge $\hat V$ if the germ of $S_i$ is contained in $\partial V$.  It follows
that $\cG$ is homeomorphic to $\partial W$. Since $f$ maps the separatrices
$S_i$
to themselves, it follows that the system of neighborhoods $W$ is also
preserved
by $f$.  Thus the structure of $\cG$ is preserved.  If $f$ preserves/reverses
orientation, $fW$ has the same/opposite orientation as $W$.  Thus we see
that we
have:

\proclaim Lemma 5.5. $\cG$ has the combinatorial structure of a simplicial
circle, and $f$ induces a simplicial homeomorphism $\hat f$ on $\cG$.   If
$f$ preserves/reverses orientation, then so does $\hat f$.

Let us note that if $f$ has the form (1.1), then
$$f(x,y)=(y,\epsilon y^d+\cdots-ax).\eqno(5.1)$$
We have $a>0$ if $f$ preserves orientation, and $a<0$ if $f$ reverses
orientation.  We conjugate by $\tau(x,y)=(\alpha x,\beta y)$ with
$\alpha,\beta\in\R$, so that
$\epsilon=\pm1$.  If $d$ is even, we require $\epsilon=+1$.  If $d$ is odd, we
define $\epsilon(f)=\epsilon$.  If $f_1$ and $f_2$ are both of odd degree
and in
the form (5.1), then $\epsilon(f_1f_2)=\epsilon(f_1)\epsilon(f_2)$.  Note that
$f^{-1}(x,y)=( a^{-1}(\epsilon x^d- y),x)$, and thus $f^{-1}$ may be put in the
form (5.1) after conjugation by $\nu(x,y)=(y,x)$ and a mapping of the form
$\tau$.  If $d$ is odd, then  $\epsilon(f^{-1})=\pm\epsilon(f)$, with the plus
sign occurring iff $f$ preserves orientation

We write $V^+=\{(x,y)\in\R^2:|y|\ge\max(R,|x|)\}$ as $V^+=V^+_1\cup
V^+_2$, with $V^+_1:=V^+\cap\{y>0\}$ and
$V^+_2:=V^+\cap\{y<0\}$.  We choose $R$ large enough that $V,V^\pm$ give a
filtration, i.e.\ $fV^+\subset V^+$ and $f(V\cup V^+)\subset V\cup V^+$.
The condition
$\epsilon(f)=1$ is equivalent to $f(V^+_1)\subset V^+_1$.  In this case, it
follows that $f(V^+_2)\subset V^+_1$ if the degree $d$ is even, and
$f(V^+_2)\subset V^+_2$ if $d$ is odd.

\proclaim Theorem 5.6.  If $f$ preserves orientation, and if $\epsilon(f)=+1$,
then all the $p_i$ are fixed.  If $\epsilon(f)=-1$, then all the $p_i$ have
period 2.

\give Proof.  Recall that for each unstably one-sided point, the separatrix
$S_i$ is unique.  Let $S_{i_1},\dots,S_{i_m}$ denote the separatrices
corresponding to $V^+_1$. That is, these are the separatrices whose germs at
infinity are contained in $V^+_1$.  Let $\cI_1\subset\cG$ denote the subgraph
whose vertices are the separatrices corresponding to $V^+_1$ and whose
edges are
the open sets $V$ with $\hat V\subset V^+_1$.  It follows that $\cI_1$ is a
proper subinterval of $\cG$. Similarly, we let $\cI_2$ denote the interval
generated by the separatrices whose germs are contained in $V^+_2$.  Note that
all of the separatrices $S_i$ are subsets of unstable manifolds.  Thus their
germs are contained in $V^+$, and so they belong to either $\cI_1$ or $\cI_2$.

If $\epsilon(f)=1$, then $fV^+_1\subset V^+_1$.  Thus $\hat f$ maps $\cI_1$ to
itself.  If $f$ preserves orientation, then $\hat f:\cI_1\to\cI_1$ is an
orientation-preserving simplicial homeomorphism.  It follows that
$\hat f$ is the identity on $\cI_1$, and so $f$ is the identity on $\cG$. Thus
$f$ maps each edge and each vertex to itself, or $f(p_i)=p_i$ for every
$i$; and this completes the proof in this case.

The remaining case is $\epsilon(f)=-1$, which implies that the degree is
odd, and thus $fV^+_1\subset V^+_2$ and $fV^+_2\subset V^+_1$.  This means that
$\hat f:\cI_1\to\cI_2$, and $\hat f:\cI_2\to\cI_1$.  Thus none of the
separatrices can be fixed.  On the other hand, $f^2$ preserves orientation and
satisfies $\epsilon(f^2)=1$.  By the argument above, these points are fixed for
$f^2$, so their periods are equal to 2. \qed

\proclaim Corollary 5.7.  All one-sided points have period 1 or 2.

\give Proof.  Without loss of generality, we may consider only
unstably one-sided points.  The mapping $f^2$ is orientation-preserving, and
$\epsilon(f^2)=+1$.  By Theorem 5.6, each one-sided point is fixed for
$f^2$.  Thus the period is 1 or 2. \qed

Let $f$ be a real, quadratic diffeomorphism of maximal entropy.  If $f$
preserves
orientation, then one of the saddle points, which we will call $p_+$, has
positive
multipliers $\lambda^u>\lambda^s>0$.  The other fixed point has negative
multipliers.

\proclaim Proposition 5.8.   Suppose $f$ is quadratic and
orientation-preserving.  Then the fixed point $p_+$ is both stably and unstably
one-sided.  No other point is one-sided.

\give Proof.  Let $p$ denote an unstably one-sided point for $f$.  By Lemma
5.6, $p$
is a fixed point.  If $p$ is the saddle point with negative multipliers, it
cannot be
one-sided.  Thus it must be the saddle point $p_+$, which has positive
multipliers.
Similarly, the stably one-sided point must also be $p_+$, so that $p_+$ is
the only
one-sided point, and it is doubly one-sided. \qed

\epsfxsize3.0in
\centerline{\epsfbox{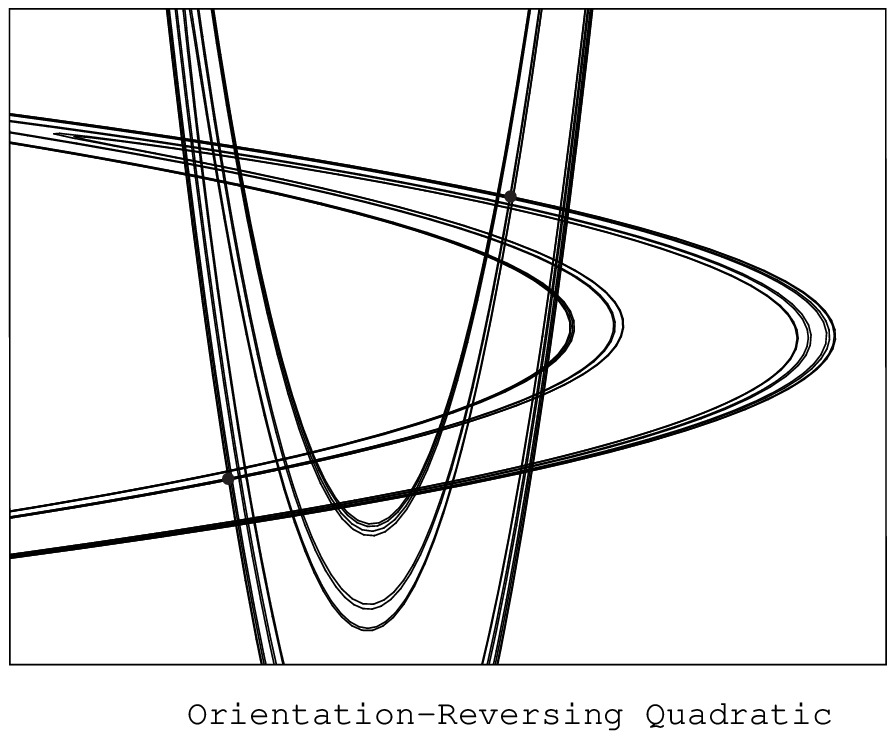}}
\centerline{Figure 5.1}

Let $f$ be a real, quadratic diffeomorphism of maximal entropy which reverses
orientation.  Then the fixed points are a pair of saddles $p_\pm$ with the
property
$\pm\lambda^u(p_\pm)>0>\pm\lambda^s(p_\pm)$.
\proclaim Proposition 5.9.  If $f$ is an orientation-reversing quadratic map,
then the one-sided points are the fixed points $p_\pm$, and and $p_\pm$ is
$u/s$-one-sided.

\give Proof.  If $p$ is a one-sided point, then the period of $p$ must be 1
or 2.  We show first that it cannot be 2.  If $q$ is a fixed point of $f$
with multipliers $\lambda^{u/s}(q)$, then $q$ is a also fixed point of $f^2$,
and it has multipliers $(\lambda^{u/s}(q))^2>0$.

Now let $h$ be an
orientation preserving map of $\R^2$.  For a saddle fixed
point $r$ of a mapping $h$, we define
$\delta(h,r)$ to be $+1$ if the multipliers of $h$ at $r$ are both positive,
and we set
$\delta(h,r)=-1$ if the multipliers are both negative.  If the degree of $h$
is even, it follows  from the Lefschetz Fixed Point Formula (cf.\
[F]) that $\sum\delta(h,r)=0$, where we sum over
all the fixed points of $h$.  Now we apply this to the map $h=f^2$, which
is orientation-preserving and degree 4.  There are 4 fixed points of $h$,
corresponding to the 2 fixed points of
$f$ and a 2-cycle for $f$.  It follows that $\{p_1,p_2\}$ is the 2-cycle
of $f$, we must have $\delta(h,p_j)=-1$ for $j=1,2$.

On the other hand, suppose that $p_1$ is one-sided with separatrix $S_1$.
Then $p_2$ is also one-sided, with separatrix $S_2$.  Now $f^2$ maps the
separatrix $S_1$ to itself, and thus the eigenvalue of $Df^2(p_1)$ must be
positive.  This means that $\delta(f^2,p_1)=+1$, which is a contradiction.
Thus the only one-sided points are the fixed points of $f$.

Since $f$ reverses orientation, each periodic point satisfies
$\lambda^u\lambda^s<0$, and thus $\lambda^s$, $\lambda^u$ cannot both be
positive.  Thus $p$ cannot be doubly one-sided.

Thus one of the fixed points must be stably one-sided with
$\lambda^s>0>\lambda^u$ for this point.  The other fixed point must be
unstably one-sided and must have $\lambda^u>0>\lambda^s$.
\qed

\proclaim Proposition 5.10.  Suppose that $f$ reverses orientation and that
$\epsilon(f)=+1$.  If $d$ is even, then at most one $u$ one-sided point can
be a
fixed point; if $d$ is odd, then at most two $u$ one-sided points can be fixed.

\give Proof.  We construct the subintervals $\cI_1$ and $\cI_2$ as in the proof
of Theorem 5.6.  If the degree of $f$ is even, then $fV^+_2\subset V^+_1$,
and thus
$\cI_2$ is empty.  The induced map $\hat f:\cI_1\to\cI_1$ is an
orientation-reversing
simplicial homeomorphism.  The simplicial map $\hat f$ can fix at most one
vertex of
$\cI_1$ (which happens when the number of vertices in
$\cI_1$ is odd).  If $d$ is odd, then $fV^+_2\subset V^+_2$, so that each
of the
restrictions $\hat f|\cI_1$ and  $\hat f|\cI_2$ can have at most one fixed
point.
\qed

\bigskip
\centerline{\bf References}
\everypar={\hangindent20pt\hangafter=1\noindent}

\bigskip
\item{[BLS]} E.\ Bedford, M.\ Lyubich, and J.\ Smillie, Polynomial
diffeomorphisms
of ${\bf C}^2$. IV: The measure of maximal entropy and laminar currents.
Invent.
math. 112, 77--125 (1993).
\item{[BS8]}  E. Bedford and J. Smillie, Polynomial diffeomorphisms
of ${\bf C}^2$. VIII: Quasi-Expansion.
\item{[BL]}  C. Bonatti and R. Langevin, {\sl Diff\'eomorphismes de Smale des
surfaces}, Ast\'erisque, 250, Soci\'ete Math\'ematique de France, 1998.
\item{[DN]} R. Devaney and Z. Nitecki, Shift automorphism in the H\'enon
mapping.
Comm.\ Math.\ Phys.\ 67, 137--148 (1979).
\item{[F]} J. Franks, {\sl Homology and Dynamical Systems} CBMS 49, AMS, 1982.
\item{[FM]} S.\ Friedland and J.\ Milnor, Dynamical properties of plane
polynomial
automorphisms. Ergodic Theory Dyn. Syst. 9, 67--99 (1989).
\item{[HO]} J.H.\ Hubbard and R.\ Oberste-Vorth, H\'enon mappings in the
complex domain II: Projective and inductive limits of polynomials, in: {\sl
Real and Complex Dynamical Systems}, B. Branner and P. Hjorth, eds. 89--132
(1995).
\item{[KKH]} I.\ Kan, H.\ Ko\c cak, J.\ Yorke
Antimonotonicity: concurrent creation and annihilation of periodic orbits. Ann.
of Math. (2) 136 , no. 2, 219--252 (1992).
\item{[M]} J. Milnor, Non-expansive H\'enon maps, Adv. Math. 69 (1988),
109--114.
\item{[NP]}  S. Newhouse and J. Palis, Hyperbolic nonwandering sets on
two-dimensional manifolds, Dyn. Systems, Peixoto ed., Salvador, 2  93--301
(1973).
\item{[O]} R. Oberste-Vorth, Complex horseshoes.  Thesis, Cornell
University 1987.
\item{[R]} C. Robinson, {\sl Dynamical Systems}, CRC Press, 1995.

\vskip.5in
\rightline{Eric Bedford}
\rightline{Indiana University}
\rightline{Bloomington, IN 47401}
\bigskip
\rightline{John Smillie}
\rightline{Cornell University}
\rightline{Ithaca, NY 14853}

\bye